\documentclass{amsart}

\usepackage{amssymb,amsthm}
\usepackage{mathrsfs}
\usepackage[all,cmtip]{xy}
\usepackage{color}
\usepackage{url}
\usepackage{enumerate}

\newcommand{\abcat}{\mathcal{A}\mathfrak{belian}\mathcal{C}\mathfrak{ategories}}
\newcommand{\abs}[1]{\left|#1\right|}
\newcommand{\addcat}[1]{\mathcal{#1}}
\newcommand{\adjunct}[4]{\xymatrix{#1 : #2 \ar@^{->}@<1pt>[r] & #3 : #4 \ar@^{->}@<1pt>[l]}}
\newcommand{\boldi}{\mathbf{i}}
\newcommand{\boldp}{\mathbf{p}}
\newcommand{\Cat}{\mathscr{C}\mathfrak{at}}
\newcommand{\cdg}[1]{C_{\text{dg}}(#1)}
\newcommand{\cdgr}[2]{C_{\text{dg}, #2}(#1)}
\newcommand{\cech}{\check{H}}
\newcommand{\coh}[1]{\mathcal{C}\mathfrak{oh}(#1)}

\newcommand{\dgcatsymb}[1]{\mathscr{#1}}
\newcommand{\dgcat}{\mathcal{DGC}\mathfrak{at}}

\newcommand{\Dgmod}{\mathscr{DG}\text{-}\mathscr{M}\mathfrak{od}}
\newcommand{\dgnat}[1]{\boldsymbol{#1}}

\newcommand{\Dperf}{\cD_{\operatorname{perf}}}
\newcommand{\Dtr}[1]{\cD^{tr}_{\!#1}}
\newcommand{\dual}[1]{#1^{\vee}}
\newcommand{\free}[1]{#1\operatorname{-free}}
\newcommand{\freeweak}[1]{#1\operatorname{-hfree}}

\newcommand{\idemcomp}[1]{#1^{\natural}}
\newcommand{\idf}[1]{\id_{#1}}

\newcommand{\Mod}{\mathscr{M}\mathfrak{od}}
\newcommand{\module}[1]{#1\operatorname{-mod}}
\newcommand{\moduledg}[1]{#1\textbf{-mod}}
\newcommand{\moduleweak}[1]{#1\!\operatorname{-hmod}}
\newcommand{\N}{\mathbb{N}}
\newcommand{\objA}[4]{\xymatrix{#1 \ar[r] & \perf{#1} \ar@^{->}@<1pt>[d]^{#3} \\ & #2 \ar@^{->}@<1pt>[u]^{#4}}}
\newcommand{\objAb}[4]{\xymatrix{#1 \ar[r] & \perf{(#1)} \ar@^{->}@<1pt>[d]^{#3} \\ & #2 \ar@^{->}@<1pt>[u]^{#4}}}
\newcommand{\objAbline}[4]{\xymatrix{#1 \ar[r] & \perf{(#1)} \ar@^{->}@<1pt>[r]^{#3} & #2 \ar@^{->}@<1pt>[l]^{#4}}}
\newcommand{\objAline}[4]{\xymatrix{#1 \ar[r] & \perf{#1} \ar@^{->}@<1pt>[r]^{#3} & #2 \ar@^{->}@<1pt>[l]^{#4}}}
\newcommand{\one}{1\!\!1}
\newcommand{\opp}[1]{{#1}^{\operatorname{op}}}
\newcommand{\perf}[1]{#1^{\operatorname{perf}}}

\newcommand{\pretr}[1]{#1^{\operatorname{pretr}}}

\newcommand{\propP}{\ensuremath{\mathbf{P}}}

\newcommand{\semifree}[1]{#1^{\operatorname{sf}}}
\newcommand{\semifreefg}[1]{#1^{\operatorname{sf}}_{\operatorname{fg}}}

\newcommand{\set}[2]{\left\{ #1\,\middle|\,#2 \right\}}
\newcommand{\setlist}[1]{\left\{ #1 \right\}}
\newcommand{\sheafF}{\mathcal{F}}
\newcommand{\sheafO}{\mathcal{O}}
\newcommand{\todo}[1]{}
\newcommand{\tododone}[1]{}
\newcommand{\tri}[1]{#1^{\operatorname{tr}}}

\newcommand{\Z}{\mathbb{Z}}

\DeclareMathOperator{\Br}{Br}
\DeclareMathOperator{\Hom}{Hom}

\DeclareMathOperator{\charac}{char}
\DeclareMathOperator{\Cone}{Cone}

\DeclareMathOperator{\id}{id}
\DeclareMathOperator{\im}{im}

\DeclareMathOperator{\lffr}{lffr}
\DeclareMathOperator{\corrmod}{(Mod)}

\DeclareMathOperator{\Pic}{Pic}
\DeclareMathOperator{\PreTr}{Pre-Tr}
\DeclareMathOperator{\Spec}{Spec}

\newcommand{\cA}{\addcat{A}}
\newcommand{\cB}{\addcat{B}}
\newcommand{\cC}{\addcat{C}}
\newcommand{\cD}{\addcat{D}}

\newcommand{\cF}{\addcat{F}}
\newcommand{\cG}{\addcat{G}}

\newcommand{\cI}{\addcat{I}}

\newcommand{\cK}{\addcat{K}}

\newcommand{\cT}{\addcat{T}}

\newcommand{\cV}{\addcat{V}}

\newcommand{\dA}{\dgcatsymb{A}}
\newcommand{\dB}{\dgcatsymb{B}}
\newcommand{\dC}{\dgcatsymb{C}}
\newcommand{\dD}{\dgcatsymb{D}}
\newcommand{\dE}{\dgcatsymb{E}}

\newcommand{\dK}{\dgcatsymb{K}}

\newcommand{\dQ}{\dgcatsymb{Q}}

\newcommand{\bC}{\mathbf{C}}
\newcommand{\bD}{\mathbf{D}}

\newcommand{\bF}{\mathbf{F}}
\newcommand{\bG}{\mathbf{G}}
\newcommand{\bH}{\mathbf{H}}
\newcommand{\bI}{\mathbf{I}}

\newcommand{\bK}{\mathbf{K}}
\newcommand{\bL}{\mathbf{L}}
\newcommand{\bM}{\mathbf{M}}

\newcommand{\dalpha}{\dgnat{\alpha}}
\newcommand{\dbeta}{\dgnat{\beta}}

\newcommand{\depsilon}{\dgnat{\epsilon}}

\newcommand{\deta}{\dgnat{\eta}}

\newcommand{\dlambda}{\dgnat{\lambda}}
\newcommand{\dmu}{\dgnat{\mu}}
\newcommand{\dnu}{\dgnat{\nu}}

\newcommand{\drho}{\dgnat{\rho}}

\numberwithin{equation}{section}

\newtheorem{theorem}[equation]{Theorem}
\newtheorem{lemma}[equation]{Lemma}
\newtheorem{proposition}[equation]{Proposition}
\newtheorem{corollary}[equation]{Corollary}
\theoremstyle{definition}
\newtheorem{definition}[equation]{Definition}

\theoremstyle{remark}
\newtheorem{remark}[equation]{Remark}
\newtheorem{example}[equation]{Example}

\begin{document} 

\title{On Differential Graded Eilenberg-Moore Construction}

\author{Umesh V. Dubey}
\address{Umesh V. Dubey, Harish-Chandra Research Institute, HBNI, Chhatnag Road, Jhusi,
Allahabad - 211 019, India}
\email{umeshdubey@hri.res.in}

\author{Vivek Mohan Mallick}
\address{Vivek Mohan Mallick, Office A-417, Main Building, Indian Institute of Science
Education and Research Pune, Dr. Homi Bhabha Road, Pune - 411 008, India}
\email{vmallick@iiserpune.ac.in}

\subjclass[2000]{Primary 14A22; Secondary 18E30 55U35}

\keywords{Differential graded categories; monads; triangulated categories;
equivariant categories; localization; twisted derived categories}

\begin{abstract}
 The Eilenberg-Moore construction for modules over a differential graded
 monad is used to study a question of Balmer regarding existence of an exact
 adjoint pair representing an exact monad. A Bousfield-like localization for
 differential graded categories is realized as a special case of this
 construction using Drinfeld quotients. As applications, we study some
 example coming from $G$-equivariant triangulated categories and twisted
 derived categories.
\end{abstract}

\maketitle


\section{Introduction}

\tododone{Mention existence of enhancements also gives
examples where we can give a description of the module category.}

Monads are ubiquitous. For example, in tensor triangulated geometry, they
were used by Balmer to characterize separated {\'e}tale morphism of
quasi-compact and quasi-separated schemes \cite{balmer:etalneth}. Kleisli
\cite{kleisli:monadadp}, and Eilenberg and Moore
\cite{eilenbergmoore:adjfunctrip} proved that any monad is a composition of
an adjoint pair of functors. In one of his papers, Balmer
\cite{balmer:septrcat} asked the question whether a monad on a triangulated
categories can be written as a composition of an adjoint pair of exact
functors.  This question is difficult to answer while staying in the world
of tensor triangulated categories. This paper gives a partial answer to this
question when the triangulated category has a suitable enhancement (see
theorem \ref{prp:factexad}). 

In the process, we do the Eilenberg-Moore construction over DG categories.
The naive generalization of the definition of a monad is too restrictive for
applications. That was the motivation for defining weak monads (definition
\ref{dfn:wkmonads}).  

One knows that Bousfield localization functors are monads by definition. We
show that weak Bousfield localization functors correspond to Drinfeld
quotients in proposition \ref{prp:bouslocw}.

As an application we have a reinterpretation of a construction done by
Sosna \cite{sosna:lintrcat} and Elagin \cite{elagin:onequivtrcat} in terms
of monads in section \ref{sec:gequivdc}. The interpretation of Elagin's result
in light of separable monads was arrived at independently in Chen
\cite{chen:sepfuncmonad}.

This paper is organized as follows. We begin by recalling some known results
in section \ref{sec:prelimns}. This section is divided into three parts: a
reveiw of Balmer's theory in \ref{ssc:balmertr}; an overview of some facts
on DG categories in \ref{ssc:dgcategory}; and finally a short subsection
dealing with monads on enriched categories in \ref{ssc:monadsenr}. Section
\ref{sec:fcttrmod} gives some conditions on exitence of triangulation of
order $N$ on modules over a monad; given that such a triangulation exists
after going to the idempotent completion.

We introduce the concept of DG monads on DG categories in section
\ref{sec:dgmonads}. Initially, in subsection \ref{ssc:monaddgc}, we study
the DG monads as a special case of monads on enriched categories. However
the definition of modules over DG monads is slightly different from that
coming from the theory of enriched categories. The slight modification
helped us get the required result. We end that subsection with a few
examples.  Subsection \ref{ssc:weakmond} weakens the definition of DG
monads. We study some properties analogues to classical monads in this set
up.

As a first application, we study a DG version of Bousfield localization in
section \ref{sec:bousfield}. The second application, in section
\ref{sec:gequivdc}, is triangulation on $G$-equivariant categories of
triangulated categories with an action of a finite group $G$. To make the
proof a bit more readable, we did a general construction, generalizing
Elagin's construction \cite{elagin:onequivtrcat} in subsection
\ref{ssc:gequivgen}. Elagin's construction is interpreted in our set up in
subsection \ref{ssc:elagincst}. We end with some applications, where we
compare $G$-equivariant derived categories with various other categories.

Section \ref{sec:twistedc} deals with the derived category of twisted
sheaves. We prove that they can be interpreted as category of modules over a
monad and using that we construct an enhancement of the twisted derived
category.  We get some more applications by considering composition of two
monads in section \ref{sec:compatmnd}. This gives an enhancement of the
derived category of twisted stheaves with cohomology supported on some
closed subvariety. We end the section by applying compatibile monads to
study the interaction of Bousfield localization and $G$-equivariance for a
finite group $G$.

\subsection*{Acknowledgement}
The first author would like to thank DST INSPIRE for funding this research
along with IISc Bangalore. He also thanks CRM Barcelona, MRC, University of
Warwick where part of the research was done. He would like to thank Prof.
Paul Balmer for his interest in the work and comments. The second author
would like to thank IISER Pune for providing him with a great ambience and
infrastructure, from where he contributed to this work. We thank Prof. Chen
for informing us about his paper \cite{chen:sepfuncmonad}.

\section{Conventions} \label{sec:convenns}
Throughout this paper, we fix the conventions mentioned in this section.  We
shall be dealing with both DG and triangulated categories. Unless otherwise
stated, we shall denote triangulated categories using the (mathcal) symbols:
$\cA$, $\cB$, $\cC$ and $\cD$. In certain parts of the text, we shall use
these same symbols to denote additive categories. $\dA\!$, $\dB$, $\dC$ and
$\dD$ (mathfrak symbols) will denote DG categories. Monads, DG functors
between DG categories and (weak) DG natural transformations between DG
functors will be denoted by bold fonts (mathbf), e.g. $\bF$, $\bG$,
$\dalpha$, $\dbeta$, etc.

\section{Preliminaries} \label{sec:prelimns}

\subsection{Balmer's theory of monads on $n$-triangulated categories}
\label{ssc:balmertr}

\tododone{Monads on triangulated categories, Separable monads, Eilenberg-Moore,
universal property, criteria for M-mod to be n-triangulated, idempotent
completion.}

In this section, all references follow the facts stated here.

\begin{definition}
  A \emph{suspended category} is an additive category $\cC$ admitting an
  auto-equivalence $T$. As in Balmer \cite[Definition 1.1]{balmer:septrcat},
  we shall assume that $T$ is an isomorphism and $T^{-1} T = T T^{-1} =
  \id_{\cC}$, the identity functor on $\cC$.

  We recall the construction of an $n$-triangulated categories and refer to
  \cite[section 5]{balmer:septrcat} for a complete description.

  For $n \geq 2$, an \emph{$n$-triangle} in a suspended category is a
  collection of objects $a_{i, j}$ with $(i, j) \in \Z \times \Z$ and a
  collection of morphisms $h_{i, j} : a_{i, j} \to a_{i, j+1}$ and $v_{i, j}
  : a_{i, j} \to a_{i+1, j}$ such that: \textit{(a)} the non-zero $a_{i,
  j}$'s are concentrated on the strip $1 \leq \abs{i - j} \leq n$;
  \textit{(b)} $a_{i, j+n+1} = Ta_{j, i}$, $h_{i, j+n+1} = Th_{j, i}$ and
  $v_{i, j+n+1} = Tv_{j, i}$; and \textit{(c)} $v_{i, j+1} \circ h_{i, j} =
  h_{i+1, j} \circ v_{i, j}$. The collection $a_{0, 1} \to a_{0, 2} \to
  \dotsb \to a_{0, n}$ is called the \emph{base of the $n$-triangle}.

  Let $N \geq 2$. An \emph{triangulation of order $N$} on a suspended
  category $\cC$ is a collection of \emph{distinguished $n$-triangles} for
  $n \leq N$, which satisfy axioms similar to all the axioms of a
  triangulated category, except for the octahedral axiom. See
  \cite[definition 5.9]{balmer:septrcat} for details.

  A \emph{functor} between categories with triangulation of order $N$ is
  \emph{exact up to order $N$} if it commutes with suspension and preserves
  distinguished $N$-triangles. By an exact functor between two categories
  with triangulation of order $N$, we shall mean a functor which is exact up
  to order $N$.

  A category with triangulations of infinite order is a collection of
  distinguished $n$-triangles for all $n \in \N$.  Such categories are
  called \emph{$\infty$-triangulated}.
\end{definition}

We need the following remark from Balmer \cite[Remark
5.15]{balmer:septrcat}.

\begin{remark}
  Every $a_{i, j}$ in a distinguished $n$-triangle $\Theta$ is a cone of the
  morphism $a_{0, i} \to a_{0, j}$, which is a composition of morphisms in
  the base of $\Theta$.
  \label{rmk:objncone}
\end{remark}

\begin{remark}
  If $\cA$ is an algebraic \cite[page 389]{schwede:algtoptrcat}
  \tododone{reference, Schwede} triangulated category, then it is the homotopy
  category of a stable model category. By Balmer \cite[remark
  5.12]{balmer:septrcat} such categories are $\infty$-triangulated.
  \label{rmk:algtrcat}
\end{remark}

\tododone{Define idempotent completion}

Now we briefly review idempotent completions. For definitions, we refer to
\cite{bs:idemcompl}.
\begin{definition}
  Suppose $\cK$ is an additive category. A morphism $e : A \to A$ is said to
  be \emph{idempotent} if $e^2 := e \circ e = e$. $\cK$ is said to be
  \emph{idempotent complete} or \emph{Karoubian} if every idempotent $e : A
  \to A$ arises from a splitting $A = \ker(e) \oplus \im(e)$. For a category
  $\cC$, its \emph{idempotent completion} $\idemcomp{\cC}$ is an idempotent
  complete category such that
  \begin{itemize}
  \item there is a fully faithful additive functor $\iota : \cC \to
    \idemcomp{\cC}$ and
  \item every additive functor $F : \cC \to \cK$ to an idempotent complete
    category factors uniquely via $\idemcomp{F} : \idemcomp{\cC} \to \cK$;
    i.e.\ $F = \idemcomp{F} \circ \iota$.
  \end{itemize}
  \label{dfn:idmcmplt} 
\end{definition}

It is well known that any additive category admits an idempotent completion.
Balmar and Schlichting \cite[Theorem 1.5]{bs:idemcompl},
proved that the idempotent completion of a triangulated category has a
unique triangulated structure such that the inclusion functor $\iota$ is
exact. We need the following generalization of that theorem.

\begin{theorem}
  If $\dC$ is a category with a triangulation of order $N$, then the
  idempotent completion $\idemcomp{\dC}$ is also a category with
  triangulation of order $N$. Moreover the inclusion functor $\iota :
  \dC \to \idemcomp{\dC}$ is exact up to order $N$.
  \label{thm:idcmpnt}
\end{theorem}

\begin{proof}
  The arguments in \cite{bs:idemcompl} can be modified to get the above
  result.
\end{proof}

\begin{corollary}
  A candidate for $\idemcomp{\cC}$ is the category whose objects are pairs
  $(A, e)$ where $A$ is an object of $\cC$ and $e : A \to A$ is an
  idempotent morphism; and whose morphissm $\varphi : (A, e) \to (B, f)$ are
  morphisms $\varphi : A \to B$ in $\cC$ such that $\varphi \circ e = f
  \circ \varphi$.
  \label{cor:descidcm}
\end{corollary}
\begin{proof}
  Follows from the proof of the theorem \ref{thm:idcmpnt}.
\end{proof}

\tododone{Define monads on $n$-triangulated categories; exact monads}

The following definition is from \cite[Definition 2.1]{balmer:septrcat}.
\begin{definition}
  A monad $(M, \mu, \eta)$ on a category $\cC$ with a triangulation of
  order $N$ is a functor $M : \cC \to \cC$ which is exact up to order $N$
  and $\mu : M \circ M \to M$ and $\eta : \idf{\cC} \to M$ are natural
  transformations such that the following diagrams commute.
  \begin{equation}
    \xymatrix{
      M \circ M \circ M \ar[r]^-{M \mu} \ar[d]_{\mu M} & M \circ M
      \ar[d]^{\mu}
      & &
      M \ar[r]^-{M \eta} \ar@{=}[dr] & M \circ M \ar[d]^{\mu} & M
      \ar[l]_-{\eta M}  \ar@{=}[dl] \\
      M \circ M \ar[r]_-{\mu} & M
      & &
      & M & 
    }
    \label{eqn:monadcmp}
  \end{equation}
\end{definition}

Given a monad $(M, \mu, \eta)$ on a category $\cC$, one can construct a
category of modules over $M$ as follows \cite[Definition
2.4]{balmer:septrcat}.

\begin{definition}
  A (left) Eilenberg-Moore category of (left) $M$ modules is the category
  whose
  \begin{itemize}
    \item objects are pairs $(x, \lambda)$ where $x$ is an object of
      $\cC$ and $\lambda : Mx \to x$ is a morphism in $\cC$ such that the
      following diagrams commute.
      \begin{equation}
        \xymatrix{
          M^2 x \ar[r]^{M \lambda} \ar[d]_{\mu} & Mx \ar[d]^{\lambda}
          & &
          x \ar[r]^{\eta_x} \ar@{=}[dr] & Mx \ar[d]^{\lambda}
          \\
          Mx \ar[r]^{\lambda} & x
          & &
          & x
        }
        \label{eqn:Mmodules}
      \end{equation}
    \item and morphisms $\varphi : (x, \lambda) \to (y, \tau)$ are morphisms
      $\varphi : x \to y$ in $\cC$ such that the following diagram commutes.
      \begin{equation}
        \xymatrix{
          Mx \ar[r]^{M \varphi} \ar[d]_{\lambda} & My \ar[d]^{\tau} \\
          x \ar[r]^{\varphi} & y
        }
        \label{eqn:morphmod}
      \end{equation}
	  The category of $M$ modules are denoted by $\module{M}$.
	\item We have a functor $F_M : \cC \to \module{M}$ defined by
	  $F_M(x) = (Mx, \mu_x)$ on objects $x$ of $\cC$ and for morphisms
	  $\varphi : x \to y$ in $\cC$, $F_M(\varphi) = M \varphi$ gives a
	  morphism between the corresponding modules.
	\item $F_M$ has a right adjoint $G_M$, which is the forgetful functor
	  defined by $G_M(x, \lambda) = x$.
  \end{itemize}
  \label{dfn:emcons}
\end{definition}

We remind ourselves of separable monads \cite[Definition
3.5]{balmer:septrcat}.

\begin{definition}
  A monad $(M, \mu, \eta)$ on a category $\cC$ is said to be
  \emph{separable} if the multiplicaiton map $\mu : M^2 \to M$ admits a
  section $\sigma : M \to M^2$, i.e.\ $\mu \circ \sigma = \idf{M}$ such that
  \begin{equation*}
    M \mu \circ \sigma M = \mu M \circ M \sigma = \sigma \circ \mu.
  \end{equation*}

  Furthermore, if $\cC$ is suspended and $M$ is stable, we call $M$ to be
  \emph{stably separable}.
  \label{dfn:sepmonad}
\end{definition}

\begin{theorem}[Balmer]
  Suppose
  \begin{enumerate}
    \item $\cC$ is an idempotent complete category with triangulation of
      order $N \geq 2$; and
    \item $(M, \mu, \eta)$ is a stably separable monad on $\cC$ such that
      $M$ is exact up to order $N$.
  \end{enumerate}
  Then we have the following
  \begin{enumerate}
    \item $\idemcomp{(\free{M})} \cong \module{M}$.
    \item $\module{M}$ admits a triangulation of order $N$. For $n \leq N$,
      an $n$-triangle $\Theta$ is distinguished in $\module{M}$ if and only
      if its image under the forgetful functor $G_M(\Theta)$ is distinguished
      in $\cC$.
    \item $F_M$ and $G_M$ are exact up to order $N$.
  \item Suppose $\adjunct{F}{\cC}{\cD}{G}$ be an adjunction such that
    $\cD$ is idempotent complete, $G$ is stably separable and
    $G \circ F = M$. Then $K$ in the following diagram is an equivalence
    and $L$ becomes an equivalence after idempotent completion.
    \begin{equation*}
    \xymatrix{
      &
      \cC \ar@<1pt>@{-^{>}}[dl] \ar@<1pt>@{-^{>}}[d] \ar@<1pt>@{-^{>}}[dr]&
      \\
      \free{M} \ar@<1pt>@{-^{>}}[ur] \ar[r]^-{L} &
      \cD \ar@<1pt>@{-^{>}}[u] \ar[r]^-{K} &
      \module{M} \ar@<1pt>@{-^{>}}[ul]
    }
    \end{equation*}
  \end{enumerate}
  \label{thm:Mmodntrg}
\end{theorem}
The above theorem is exactly parts (a) - (d) of \cite[Main Theorem
5.17]{balmer:septrcat}. \tododone{Insert part (d) only if it is relevant to
this paper.}

\subsection{DG categories} \label{ssc:dgcategory}

\tododone{Definition, functors, 2-cat of DG categories, $H^0$, Pre-Tr, Tr,
  (strong) enhancement, dg-lifts, pre-exact functors, quasi-equivalence,
  semi-free, perfect DG categories as an analogue of idempotent completion
  for triangulated categories.\\ (Keller, Drinfeld, Orlov-Lunts,
Bondal-Orlov)}

We use the following definition of a DG category. For further details we
refer to \cite[\S 1]{bondalkapranov:enhancedtrcat} and \cite[section
2]{drinfeld:dgquotdgcat}.

\begin{definition}
  Suppose $k$ is a commutative ring. A \emph{DG category}, or a
  \emph{differential graded $k$-linear category}, is a category such that
  \begin{itemize}
    \item for objects $A$ and $B$ in $\dC$, $\Hom_{\dC}(A, B)$ has
      \begin{itemize}
        \item a structure of a $\Z$-graded $k$ module,
        \item a differential $d : \Hom_{\dC}(A, B) \to \Hom_{\dC}(A, B)$ of
          degree $1$ (with $d^2 = 0$);
      \end{itemize}
    \item for objects $A$, $B$ and $C$ of $\dC$, the composition morphism
      \begin{equation*}
        \Hom_{\dC}(B, C) \otimes_k \Hom_{\dC}(A, B) \to \Hom_{\dC}(A, C)
      \end{equation*}
      is a morphism of complexes; and
    \item $d(\idf{A}) = 0$ for every object $A$.
  \end{itemize}
  Given objects $A$ and $B$, the group of all morphisms of degree $r$
  between $A$ and $B$ is denoted by $\Hom_{\dC}^r(A, B)$.
  \label{def:dgcatgry}
\end{definition}

\begin{definition}
  A \emph{DG functor} $\bF$ between DG categories $\dC$ and $\dD$ is a
  $k$-linear functor such that
  \begin{itemize}
    \item $\deg \bF f = \deg f = r$ for all morphisms $f \in \Hom^r_{\dC}(A,
      B)$ in $\dC$, and
    \item $\bF (d f) = d (\bF f)$ for all morphisms $f$ in $\dC$.
  \end{itemize}
\end{definition}

\begin{definition}
  Suppose $\dC$ and $\dD$ are two DG categories and $\bF$ and $\bG$ be two
  DG functors between them. Given a natural transformation $\dalpha : \bF
  \to \bG$, $d \dalpha$ is the natural transformation such that for every
  object $A$ in $\dC$, $(d \dalpha) A = d (\dalpha(A))$. A \emph{DG natural
  transformation} is a natural transformation $\dalpha : \bF \to \bG$ such
  that
  \begin{itemize}
    \item for each object $A$ in $\dC$, $\dalpha(A) : \bF A \to \bG A$ is of
      degree $0$, and
    \item $d \dalpha = 0$.
  \end{itemize}
  \label{dfn:dgnattrn}
\end{definition}

The fact that $\Hom_{\dC}(A, B)$ is a complex for objects $A$ and $B$ in a
DG category $\dC$, allows us to make the following definition.

\begin{definition}
  Given a DG category $\dC$, 
  $H^0(\dC)$ is the $k$-linear category whose objects are the same as those
  of $\dC$, and whose morphisms are
  \begin{equation*}
    \Hom_{H^0(\dC)}(A, B) = H^0(\Hom_{\dC}(A, B))
  \end{equation*}
  for any pair of objects in $\dC$. Given a DG functor $\bF : \dC \to \dD$
  between two DG categories, $H^0(\bF) : H^0(\dC) \to H^0(\dD)$ is the
  functor defined by
  \begin{eqnarray*}
    H^0(\bF)(A) &=& \bF(A) \text{ for all objects } A \text{ in } \dC, \\
    H^0(\bF)([f]) &=& H^0(\bF f)
  \end{eqnarray*}
  where $[f]$ is a morphism in $\Hom_{H^0(\dC)}(A, B)$ and $f$ is a
  representative in $\Hom_{\dC}(A, B)$ such that $df = 0$ and $[f] = f \mod
  d \Hom^{-1}_{\dC}(A, B)$.
  For a DG natural transformation $\dalpha : \bF \to \bG$, one defines
  $H^0(\dalpha)$ to be the natural transformation from $H^0(\bF)$ to
  $H^0(\bG)$ such that for every object $A$ in $\dC$,
  \begin{gather*}
    H^0(\dalpha)(A) = \dalpha(A) \mod d\Hom^{-1}_{\dC}(\bF(A), \bG(A)) \in
    H^0(\Hom_{\dD}(\bF(A), \bG(A))) \\
    = \Hom_{H^0(\dD)}(H^0(\bF(A)),
    H^0(\bG(A))).
  \end{gather*}
\end{definition}

\begin{definition}
  Small DG categories along with DG functors and DG natural transformations
  form a 2-category in the sense of \cite[Chapter XII]{maclane:catworkmath}.
  We denote this category by $\dgcat$.
  \label{dfn:catdgcat}
\end{definition}

\begin{lemma}
  $H^0$ defines a strict 2-functor from $\dgcat$ to the 2-category of small
  $k$-linear categories.
  \label{lem:h0strict}
\end{lemma}

\begin{proof}
  It is a routine checking of the defintion of a strict 2-functor.
\end{proof}

The following definitions are by Bondal and Kapranov
\cite{bondalkapranov:enhancedtrcat}.

\begin{definition}
  Consider a collection of objects $\set{A_i}{i \in \Z}$ and a collection of
  morphisms $\set{q_{i, j} : A_i \to A_j}{(i, j) \in \Z \times \Z}$ in some
  DG category $\dA$. Such a collection of objects and morphims is called a
  \emph{twisted complex} over $\dA$ if
  \begin{itemize}
    \item all but finitely many of $A_i$'s are $0$,
    \item $\deg q_{i, j} = i - j + 1$, and
    \item $d q_{i, j} + \sum_{k} q_{k, j} \circ q_{i, k} = 0$.
  \end{itemize}
  Let $T = (A_i, q_{i, j})$ and $U = (B_i, r_{i, j})$ be two twisted
  complexes. A morphism $f : T \to U$ of degree $k$ is a collection of
  morphisms $f_{i, j} : A_i \to B_j$ with $\deg f_{i, j} = i - j + k$.  The
  twisted complexes over a DG category form a DG category with respect to a
  suitably defined differential \cite[Definition
  1]{bondalkapranov:enhancedtrcat}.  Let $\overline{\dA}$ be the DG category
  obtained from $\dA$ by adjoining finite formal direct sums of objects.
  Define $\pretr{\dA}$ to be the category of twisted complexes over
  $\overline{\dA}$. For two objects $T = (A_i, q_{i, j})$ and $U = (B_i,
  r_{i, j})$ in $\pretr{\dA}$ and a morphism of degree $0$, $f = (f_{i, j})
  : T \to U$, such that $df = 0$, one can define the \emph{cone},
  \begin{equation*}
    \Cone(f) = (C_i, s_{i, j}) \quad \text{where } C_i = A_i + B_{i-1},\ 
    s_{i, j} =
    \begin{pmatrix}
      q_{i, j} & f_{i, j} \\
      0 & r_{i, j}
    \end{pmatrix}.
  \end{equation*}
  $\pretr{\dA}$ is called the \emph{pretriangulated hull} of $\dA$
  (see \cite[section 2.3]{orlov:gluingdgcat}). There is a canonical functor
  $\dA \hookrightarrow \pretr{\dA}$.
  \label{dfn:twistdcp}
\end{definition}

\begin{definition}
  Define $\tri{\dA} = H^0(\pretr{\dA})$. Let us define the
  \emph{distinguished triangles} in $\tri{\dA}$ to be the triangles induced
  by the natural morphisms
  \begin{equation*}
    A \xrightarrow{f} B \rightarrow \Cone{f} \xrightarrow{+}
  \end{equation*}
  in $\pretr{\dA}$, for every $f$ of degree $0$ satisfying $df = 0$; and all
  triangles isomorphic to these.
\end{definition}

\begin{proposition}
  $\tri{\dA}$ along with the distinguished triangles defined above forms a
  triangulated category.
  \label{prp:h0pretrtr}
\end{proposition}
\begin{proof}
  This is \cite[Proposition 1]{bondalkapranov:enhancedtrcat}.
\end{proof}

\begin{definition}
  A DG functor $\bF : \dA \to \dB$ is said to be a \emph{quasi-equivalence}
  if for every pair of objects $A$ and $A'$ in $\dA$, the morphism $\bF_{A,
  A'} : \Hom_{\dA}(A, A') \to \Hom_{\dB}(\bF A, \bF A')$ is a
  quasi-isomorphism and if the induced functor $H^0(\bF) : H^0(\dA) \to
  H^0(\dB)$ is an equivalence. Two DG categories $\dA$ and $\dB$ are said to
  be \emph{quasi-equivalent} if there exists DG categories $\dC_1, \dotsc,
  \dC_n$ and a chain of quasi-equivalences:
  
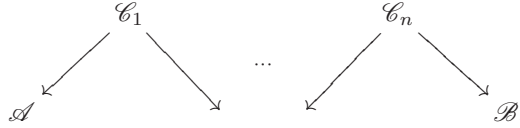
\begin{equation*}
    \xymatrix{
      & \dC_1 \ar[dr] \ar[ld] & \ar@{}[dr]|{\cdots} &         & \dC_n
      \ar[dr] \ar[ld] & \\
      \dA & & & & & \dB.
    }
  \end{equation*}
\end{definition}

\begin{definition}
  Two objects $c$ and $d$ in a DG category $\dC$ are said to be
  \emph{homotopy equivalent} if there are morphisms $f : c \to d$ and
  $g : d \to c$ in $\dC$ such that $g \circ f - \idf{c} = d \theta$ and
  $f \circ g - \idf{d} = d \tau$ for morphisms $\theta$ and $\tau$ in $\dC$.
  \label{def:homequiv}
\end{definition}

\begin{definition}
  A DG category $\dA$ is said to be \emph{pretriangulated} if the canonical
  functor $\dA \to \pretr{\dA}$ is a quasi-equivalence. 
  A \emph{strongly pretriangulated} DG category is a DG category $\dA$ for
  which the canonical functor $\dA \to \pretr{\dA}$ is a DG-equivalence.
\end{definition}

\begin{remark}
  An equivalent condition for a non-empty DG category $\dA$ to be
  pretriangulated is that
  \begin{itemize}
    \item for any object $A$ in $\dA$, and for any $n \in \Z$, the object
      $A[n]$ in $\pretr{\dA}$ is homotopy equivalent to some object in
      $\dA$, and
    \item for any closed degree zero morphism $f$ in $\dA$, $\Cone(f)$ in
      $\pretr{\dA}$ is homotopy equivalent to some object in $\dA$.
  \end{itemize}

  If we replace ``homotopy equivalent'' with ``DG isomorphic'' in the above
  conditions we get an equivalent description of strongly pretriangulated DG
  categories. See \cite[section 2.4]{drinfeld:dgquotdgcat}.
  \label{rmk:altdscpt}
\end{remark}

We need the following result about functoriality of cones.

\begin{proposition}
  Suppose $\dC$ is a strongly pretriangulated DG category. Then any
  commutative square
  \begin{equation*}
    \xymatrix{
      A \ar[r]^{f} \ar[d]_{g} &
      B \ar[d]^{h} \\
      C \ar[r]^{k} &
      D
  }
  \end{equation*}
  with $f$ and $k$ being degree zero closed morphisms, one has a canonical
  way to fill the thrid arrow below to make all the squares commutative.
  \begin{equation*}
  \xymatrix@R=0.7em{
    A \ar[rr]^{f} \ar[ddd]_{g} &
    &
    B \ar[dl]  \ar[ddd]^{h} \\
    &
    \Cone(f) \ar[lu]^-{+} \ar@{.>}[d] &
    \\
    &
    \Cone(k) \ar[ld]_-{+} &
    \\
    C \ar[rr]_{k} &
    &
    D \ar[lu] \\
  }
  \end{equation*}
  \label{prp:funccone}
\end{proposition}
This follows from construction of a cone \cite[definition
2]{bondalkapranov:enhancedtrcat}. For related discussion, also see \cite[\S
5]{toen:lectdgcat}.

Most of this work deals with triangulated categories which are realized as
homotopy category of some pretriangulated category. We recall the following
definition.

\begin{definition}
  \emph{Enhancement} of a triangulated category $\cA$ is a pair $(\dA,
  \varepsilon)$ where $\dA$ is a pretriangulated DG category and
  $\varepsilon : H^0(\dA) \to \cA$ is an equivalence of triangulated
  categories. (See \cite{luntsorlov:uniqueenhancement}.)
\end{definition}

In one of the applications we need the DG category $\perf{\dC}$ associated
to a DG category $\dC$. We define it in the following sequence of
definitions. See \cite[section 2.2]{orlov:gluingdgcat} for further details.

\begin{definition}
  A \emph{right DG $\dA$-module} associated to a small DG category
  $\dA$ is a DG functor $M : \opp{\dA} \to \Dgmod(k)$ where
  $\Dgmod(k)$ is the DG category of DG $k$-modules. The category of
  DG $\dA$-modules is denoted by $\Dgmod(\dA)$.
\end{definition}

\begin{remark}
  $\Dgmod(\dA)$ is a DG category (see \cite[section 1.2]{keller:derdgcat}).
  There exists a fully faithful Yoneda DG functor $h^{\bullet} : \dA \to
  \Dgmod(\dA)$ defined by
  \begin{equation*}
    h^{\bullet}(Y) = h^Y \text{, where } h^Y(X) = \Hom_{\dA}(X, Y)
  \end{equation*}
  for objects $X$ and $Y$ in $\dA$.
\end{remark}

\begin{definition}
  A DG $\dA$-module $M$ is \emph{representable} if it is of the form $h^Y$
  for some object $Y$ in $\dA$. A DG $\dA$-module is \emph{free} if it is of
  the form $\oplus_{i = 1}^n h^{Y_i}[s_i]$ for some $n \in \N$ and for some
  shifts $s_i \in \Z$. A DG $\dA$-module $S$ is said to be \emph{semi-free}
  if it admits a filtration $0 = T_0 \hookrightarrow T_1 \hookrightarrow
  \dotsb \hookrightarrow T_k \hookrightarrow \dotsb S$ such that $T_i /
  T_{i-1}$ is free for $1 \leq i \leq k$. The full DG subcategory of
  $\Dgmod(\dA)$ consisting of all semi-free DG modules is denoted by
  $\semifree{\dA}$. A semi-free module $S$ is said to be \emph{finitely
  generated}, if $T_k = S$ for some finite $k \in \N$. The full subcategory
  of $\Dgmod(\dA)$ consisting of finitely generated semi-free modules is
  denoted by $\semifreefg{\dA}$. The \emph{perfect DG module} is an object
  of $\semifree{\dA}$ which is homotopy equivalent to a direct summand of a
  finitely generated semi-free DG module. The full subcategory of all
  perfect DG $\dA$-modules is denoted by $\perf{\dA}$.
  \label{dfn:perfdgct}
\end{definition}

\begin{remark}
  $\perf{\dA}$ is a pretriangulated DG category and contains $\pretr{\dA}$.
  Also $H^0(\perf{\dA}) \cong \idemcomp{\big(H^0(\pretr{\dA})\big)} =
  \idemcomp{\big(\tri{\dA}\big)}$. For details, see \cite[section
  2.3]{orlov:gluingdgcat}. \tododone{find a better reference.}
\end{remark}

\begin{lemma}
  The functor $\PreTr : \dgcat \to \dgcat$ defined by $\PreTr(\dA) =
  \pretr{\dA}$ is a strict 2-functor.
  \label{lem:pretr2fn}
\end{lemma}

\begin{proof}
  Follows easily from definitions.
\end{proof}

We record the following definition \cite[\S 3, definition
2]{bondalkapranov:enhancedtrcat} for later use. For a pretriangulated
category $\dA$, proposition 1 in \cite[\S
3]{bondalkapranov:enhancedtrcat} asserts the existence of a equivalence
$T_{\dA} : \pretr{\dA} = \PreTr{\dA} \to \dA$, called the
\emph{convolution}.

\begin{definition}
  A \emph{pre-exact} functor $F : \dA \to \dB$ between two pretriangulated
  categories is a DG functor such that the following diagram commutes
  \begin{equation*}
    \xymatrix{
      \pretr{\dA} \ar[r]^{\pretr{F}} \ar[d]_{T_{\dA}} &
      \pretr{\dB} \ar[d]^{T_{\dB}} \\
      \dA \ar[r]^{F} & \dB.
    }
  \end{equation*}
\end{definition}
The only result we need about pre-exact functors is the following easy fact.
\begin{lemma}
  If $F : \dA \to \dB$ is a pre-exact functor between pretriangulated
  categories, then
  \begin{equation*}
    H^0(F) : H^0(\dA) \to H^0(\dB)
  \end{equation*}
  is an exact functor between triangulated categories.
  \label{lem:prexisex}
\end{lemma}

\subsection{Monads in $\cV$-categories} \label{ssc:monadsenr}
\tododone{(Kelly)}

DG-categories can be viewed as categories enriched by $\Dgmod(k)$. Since
monads in enriched categories are already studied, we recall enriched
categories here.  The references are \cite{kelly:enrichcatth}
\cite{eilenbergkelly:closedcat} \cite{kock:enrichedmonads}
\cite{kock:clcatgenenrmonads}.

Let $\cV$ be a symmetric monoidal category, consisting of a functor $\otimes
: \cV \times \cV \to \cV$, and an object $\one_{\cV}$ which acts as an
identity for $\otimes$ satisfying certain compatibility conditions as in
\cite[Part II, Chapter 1, page 375]{levine:mixedmotives}.

\begin{definition}
  A $\cV$-category $\cC$ consists of the following data:
  \begin{enumerate}
    \item a collection of objects of $\cC$;
    \item for each pair of objects $A$ and $B$ in $\cC$, an object
      $\Hom_{\cC}(A, B)$ of $\cV$;
    \item for each triple $A$, $B$ and $C$ of objects in $\cC$, a
      ``composition'' morphism
      \begin{equation*}
        \circ_{A, B, C} : \Hom_{\cC}(B, C) \otimes \Hom_{\cC}(A, B) \to
        \Hom_{\cC}(A, C)
      \end{equation*}
      in $\cV$; and
    \item for each object $A$ in $\cC$ a morphism $\id_A : \one_{\cC} \to
      \Hom_{\cC}(A, A)$
  \end{enumerate}
  which satisfy axioms corresponding to associativity and compatibility with
  the identity morphism.
  
  A $\cV$-functor $F : \cA \to \cB$ between two $\cV$-categories is an
  association of objects together with morphisms $F(A, B) : \Hom_{\cA}(A, B)
  \to \Hom_{\cB}(F(A), F(B))$ which is compatible with composition and
  identity, for any pair of objects $A$ and $B$ in $\cA$.

  A $\cV$-natural transformation $\theta : F \to G$ between two
  $\cV$-functors $F, G : \cA \to \cB$ is a collection of maps
  $\theta(A) : \one_{\cB} \to \Hom_{\cB}(F(A), G(A))$ for each object
  $A$ in $\cA$ which is compatible with morphisms. Two natural
  transformations $\rho$ and $\theta$ can be composed to give:
  \begin{multline*}
    \rho \circ \theta (A) : \one_{\cB} \to \one_{\cB} \otimes
    \one_{\cB} \xrightarrow{\rho \otimes \theta} \Hom_{\cB}(G(A),
    H(A)) \otimes \Hom_{\cB}(F(A), G(A)) \\
    \xrightarrow{\circ_{F(A), G(A), H(A)}} \Hom_{\cB}(F(A), H(A)).
  \end{multline*}
\end{definition}

\begin{definition}
  Suppose $F : \cA \to \cB$ and $G : \cB \to \cA$ be two $\cV$-functors
  between $\cV$-categories. $F$ is said to be a left adjoint for $G$ if
  there is a $\cV$-natural isomorphism
  \begin{equation*}
    a(A, B) : \Hom_{\cB}(FA, B) \to \Hom_{\cA} (A, GB)
  \end{equation*}
  between $\cV$-functors from $\opp{\cA} \otimes_{\cV} \cB \to \cV$. For
  further details one can refer to \cite{kelly:concenrichcat}.
\end{definition}

\begin{definition}
  A $\cV$-monad $(M, \mu, \eta)$ on a $\cV$-category $\cC$ consists of an
  $\cV$-endofunctor $M : \cC \to \cC$ and two $\cV$-natural transformations
  $\mu : M \circ M \to M$ and $\eta : \id_{\cC} \to M$ such that the
  following two diagrams commute:
  \[
    \xymatrix{
      M^3 \ar[r]^{M \mu} \ar[d]_{\mu_M} & M^2 \ar[d]^{\mu} & & M \ar@{=}[dr]
      \ar[r]^{\eta_M} & M^2 \ar[d] & M \ar@{=}[dl] \ar[l]_{M\eta} \\
      M^2 \ar[r]^{\mu} & M & & & M &
    }
  \]
  \label{dfn:enrmonad}
\end{definition}

The following lemma is well known.
\begin{lemma}
  Let $\cV$ be a symmetric monoidal additive category. Let $\cC$ be a
  $\cV$-category. Then for every $\cV$-monad, $M$, there exists a
  $\cV$-category $\module{M}$, and two $\cV$-functors $F_M : \cC \to
  \module{M}$ and $G_M : \module{M} \to \cC$, such that $M = G_M \circ F_M$,
  $G_M$ is right adjoint to $F_M$ and for any other $\cV$-category $\cD$
  admitting $\cV$-functors $F : \cC \to \cD$ and $G : \cD \to \cC$, such
  that $G \circ F = M$, there is a unique $\cV$-functor $\Phi : \cD \to
  \module{M}$ as in the following diagram.  Furthermore, there is a
  $\cV$-category $\free{M}$ and a unique $\Psi : \free{M} \to \cD$ such that
  the whole diagram, given below, commutes.
  \begin{equation*}
    \xymatrix{
      & \cC \ar@^{->}@<1pt>[d]^F \ar@^{->}@<1pt>[dl]^{F_M}
      \ar@^{->}@<1pt>[dr]^{F_M} & \\
      \free{M} \ar@^{->}@<1pt>[ur]^{G_M} \ar[r]^-{\Psi} & \cD
      \ar@^{->}@<1pt>[u]^G \ar[r]^-{\Phi} & \module{M}.
      \ar@^{->}@<1pt>[ul]^{G_M}
    }
  \end{equation*}
  \label{lem:emenrmnd}
\end{lemma}

\begin{proof}
  The proof goes parallel to the Eilenberg-Moore construction
  \cite{eilenbergmoore:adjfunctrip} and \cite{kleisli:monadadp}. See for
  example, \cite{linton:monademk}.
\end{proof}

\section{A fact about triangulation in M-mod} \label{sec:fcttrmod}
\tododone{Lemma 2.2 -- Cor 2.6}

In this section the canonical functor from a category to its idempotent
completion will be denoted by $\iota$.

\begin{lemma}
  Given an exact functor $F : \cC \to \cD$ between two categories with
  triangulation of order $N$ (respectively, triangulated categories), one
  can canonically extend the functor to a functor $\idemcomp{F} :
  \idemcomp{\cC} \to \idemcomp{\cD}$, between the idempotent completions of
  $\cC$ and $\cD$, such that $\iota \circ F = \idemcomp{F} \circ \iota$,
  that is, the following diagram commutes:
  \begin{equation*}
    \xymatrix{
      \cC \ar[d]_{\iota} \ar[r]^{F} & \cD \ar[d]^{\iota} \\
      \idemcomp{\cC} \ar[r]^{\idemcomp{F}} & \idemcomp{\cD}.
    }
  \end{equation*}
  \label{lem:idemcomp}
\end{lemma}

\begin{proof}
  That the categories $\idemcomp{\cC}$ and $\idemcomp{\cD}$ have
  triangulations of order $N$, follows from a straight-forward modification
  of the proof in \cite[section 1.15]{bs:idemcompl}. By
  definition (see \cite[definition 1.10]{bs:idemcompl})
  $\idemcomp{F}(a, \varepsilon) = (F(a), F(\varepsilon))$. From this
  description the commutativity of the square is evident.
\end{proof}

\begin{lemma}
  Given an exact monad $(M, \mu, \eta)$ on a category $\cC$, with triangulations
  of order $n$, (respectively, a triangulated category), there is a monad
  $(\idemcomp{M}, \idemcomp{\mu}, \idemcomp{\eta})$ on the idempotent
  completion $\idemcomp{\cC}$ of $\cC$. $\idemcomp{M}$ is compatible with
  $M$ in the sense that $\iota \circ M = \idemcomp{M} \circ \iota$.
  \label{lem:idemcmnd}
\end{lemma}

\begin{proof}
  This is an easy consequence of lemma \ref{lem:idemcomp}.
\end{proof}

The next proposition has two versions. For brevity, we state a statement
which we call ``proposition(\propP)'', which is a proposition regarding
property \propP{} of some pretriangulated categories..
\begin{definition}
  By \emph{proposition(\propP{})}, we mean the following statement:

  \begin{quote}
    ``Suppose $\cC$ and $\cD'$ are categories satisfying property \propP{} and
    $\iota : \cC \to \idemcomp{\cC}$ is an idempotent completion of $\cC$.
    Let $\adjunct{F}{\idemcomp{\cC}}{\cD'}{G}$ be an exact adjunction. By
    $\cD$, we denote the full subcategory of $\cD'$, consisting of all
    objects $d$ in $\cD'$ such that $Gd$ is isomorphic to an object in the
    image of $\iota$.  Under these conditions we have the following:
    \begin{enumerate}[(a)]
      \item $\cD$ has property \propP{},
      \item If $F' = F \circ \iota$, and $F'$ maps $\cC$ into $\cD$, there
        exists an exact functor $G' : \cD \to \cC$, which is a right adjoint
        to $F'$:
        \begin{equation*}
          \adjunct{F'}{\cC}{\cD}{G'}.
        \end{equation*}
	  \item Let $(M, \mu, \eta)$ be a monad on $\cC$ and let $(\idemcomp{M},
		\idemcomp{\mu}, \idemcomp{\eta})$ be the monad on $\idemcomp{\cC}$
		constructed in lemma \ref{lem:idemcmnd}. If we conisder $\cD' =
		\module{\idemcomp{M}}$, $F = F_{\idemcomp{M}}$ and $G =
		G_{\idemcomp{M}}$ (notation coming from definition
		\ref{dfn:emcons}), then $\cD = \module{M}$.
    \end{enumerate}
    In particular, $\module{M}$ has property \propP{} in this case.''
  \end{quote}
\end{definition}

\begin{proposition}
  Proposition(\propP{}) is true in the following two cases.
  \begin{enumerate}
    \item \propP{} means that a category is triangulated.
    \item \propP{} means that a category admits a triangulation up to order
      $N$, for $N \in \N$ or $N = \infty$.
  \end{enumerate}
  \label{prp:redtocmp}
\end{proposition}

\begin{proof}
  We prove each part separately.
  \begin{enumerate}
    \item A full subcategory of a triangulated category is triangulated if
      and only if it is closed under shifts and every morphism can be
      completed into a triangle.

      \begin{enumerate}
        \item\label{itm:pfdtrian} Suppose $d$ is an object of $\cD$. Let $Gd
          \cong \iota(c)$. Since $G(d[n]) = G(d)[n] = \iota(c)[n] =
          \iota(c[n])$, $d[n]$ is an object of $\cD$.

      Suppose $x \xrightarrow{f} y$ be a morphism in $\cD$. Since $\cD'$
      is triangulated, there exists a distinguished triangle $x \to y
      \to z \xrightarrow{+}$ in $\cD'$. Thus $Gf$ fits into a
      distinguished  triangle $Gx \xrightarrow{Gf} Gy \to Gz$ in
      $\idemcomp{\cC}$. Let $f' : x' \to y'$ be a morphism in $\cC$ such
      that $Gx \cong \iota x'$ and $Gy \cong \iota y'$, and $Gf \cong
      \iota f'$. Since $\cC$ is a triangulated category, there exists an
      object $z'$ in $\cC$ such that $x' \xrightarrow{f'} y' \to z'
      \xrightarrow{+}$ is a distinguished triangle. Consider the diagram
          \begin{equation*}
            \xymatrix{
              Gx \ar[r]^{Gf} \ar[d]_{\cong} & Gy \ar[r] \ar[d]^{\cong} & Gz
              \ar@{.>}[d] \ar[r]^{+} & \\
              \iota x' \ar[r]^{\iota f'} & \iota y' \ar[r] & \iota z'
              \ar[r]^{+} & \\
            }
          \end{equation*}
          and by definition of triangulated categories the dotted arrow
          exists.  It is also clear that it has to be an isomorphism in
          $\idemcomp{\cC}$.  Thus $z$ is an object of $\cD$.

        \item\label{itm:pfadjftr} Note that $G$ restricts to a functor $G' :
          \cD \to \cC$, by definition of $G$. Since $F$, $G$ and $\iota$ are
          exact, so are $F'$ and $G'$. $F'$ is a left adjoint to $G'$, as
          $F$ is a left adjoint to $G$.

    \item\label{itm:pfdismod} Since $\iota : \cC \to \idemcomp{\cC}$ is
      fully faithful (see definition \ref{dfn:idmcmplt} and theorem
      \ref{thm:idcmpnt}), so is the induced functor $\module{M} \to
      \module{\idemcomp{M}}$. Now if $G_{\idemcomp{M}}(x, \lambda) = x$
      is isomorphic to an object $\iota x'$ of $\cC$, Since $\iota$ is
      fully faithful, $\lambda : \idemcomp{M} x \to x$ determines a
      $\lambda' : M x' \to x'$ and $\iota(x', \lambda') \cong (x,
      \lambda)$ and hence $\cD$ is equivalent to $\module{M}$.
      \end{enumerate}
      In the set up of item \ref{itm:pfdismod} above, the fact that
      $\module{M}$ is triangulated follows from part \ref{itm:pfdtrian}.

    \item Note that the proof of \ref{itm:pfdtrian} implies that $\cD$ is
      pretriangulated. In particular, for every morphism $f : A \to B$ in
      $\cD$, if $A \xrightarrow{f} B \to C \xrightarrow{+}$ is a
      distinguished 2-triangle in $\cD'$, then $C$ is also an object of
      $\cD$.
      \begin{enumerate}
        \item By remark \ref{rmk:objncone}, given a base $a_{0, 1} \to a_{0,
          2} \to \dotso a_{0, n}$ in $\cD$, the objects $a_{i, j}$ in the
          corresponding distinguished $n$-triangle in $\cD'$ also belong to
          $\cD$. Thus, if $\cD'$ admits a triangulation up to order
          $N$, then so does $\cD$.
        \item This proof is exactly the same as that of \ref{itm:pfadjftr}.
        \item An argument similar to that of \ref{itm:pfdismod} using
          theorem \ref{thm:idcmpnt} will give us the desired result.
      \end{enumerate}
  \end{enumerate}
\end{proof}

\tododone{Insert a remark saying how the following result is Balmer's result
  along with the observation that the general case can be reduced to the
case where the categories are idempotent complete.}

We use the above proposition to slightly modify a theorem of Balmer,
recalled here as theorem \ref{thm:Mmodntrg}.

\begin{proposition}
  If $\cC$ is a category with a triangulation up to order $N$. Suppose $M$
  is a stably separable monad on $\cC$. Then $\module{M}$ has a
  triangulation up to order $N$ with respect to which $F_M$ and $G_M$ are
  exact functors.
\end{proposition}

\begin{proof}
  \tododone{refer to the definition of separable monads. Also cite a
  reference that idempotent completion of an n-triangulated category is
  n-triangulated; Balmer Schlichting + Balmer separability paper?}

  That $M$ is stably separable (definition \ref{dfn:sepmonad}) will imply
  that there is a splitting $\sigma : M \to M^2$ of $\mu : M^2 \to M$. One
  can extend $\sigma$ to a natural transformation $\idemcomp{\sigma} :
  \idemcomp{M} \to (\idemcomp{M})^2$, as (see corollary \ref{cor:descidcm})
  \begin{equation*}
    \idemcomp{\sigma}_{(x, i)} : \idemcomp{M}(x, i) = (Mx, Mi)
    \longrightarrow (\idemcomp{M})^2 (x, i) = (M^2 x, M^2 i).
  \end{equation*}
  It is clear that $\idemcomp{\sigma}$ splits $\idemcomp{\mu}$ and satisfies
  the conditions in definition \ref{dfn:sepmonad}. Thus $\idemcomp{M}$ is
  separable; and hence stably separable.

  Applying theorem \ref{thm:Mmodntrg} to $(\idemcomp{M}, \idemcomp{\mu},
  \idemcomp{\eta})$ on $\idemcomp{\cC}$, we get that
  $\module{\idemcomp{M}}$ has a triangulation up to order $N$, with respect
  to which $F_{\idemcomp{M}}$ and $G_{\idemcomp{M}}$ are exact. Now applying
  part 2 of proposition \ref{prp:redtocmp}, we conclude that
  $\module{M}$ also has a triangulation up to order $N$, with respect to
  which $F_M$ and $G_M$ are exact.
  \tododone{We need a version of the above proposition for
  $n$-trianglulated categories}
\end{proof}

\section{DG monads} \label{sec:dgmonads}

\subsection{Monads on DG categories} \label{ssc:monaddgc}
\tododone{2.8 -- 2.14}

\tododone{Do everything in DG level. Weak stuff go to the next subsection.}

\tododone{Construction, universal property, Criterion that the DG-modules enhance
M-mod (for use in Bousfield localization type situation). Lemma 2.9, Define
$\cD^{tr}$, Propn 2.13.}

Note that DG categories are nothing but $\Dgmod(k)$ enriched categories.
This allows us to define
\begin{definition}
  A DG monad is a monad on a $\Dgmod(k)$ enriched category defined as in
  definition \ref{dfn:enrmonad}.
\end{definition}

\begin{definition}
  Let $\cC$ be a triangulated category. Suppose $(\dC, \epsilon)$ is an
  enhancement of $\cC$. Let $(M, \mu, \eta)$ be a monad on $\cC$. A
  \emph{lift of $M$ to $\dC$} is a pre-exact DG monad $(\bM, \dmu, \deta)$
  on $\dC$ such that the following diagram of functors
  \begin{equation*}
    \xymatrix{
      H^0(\dC) \ar[r]^{H^0(\bM)} \ar[d]_{\epsilon} & H^0(\dC)
      \ar[d]^{\epsilon} \\
      \cC \ar[r]^{M} & \cC,
    }
  \end{equation*}
  and the following diagrams of natural transformations commute.
  \begin{equation*}
    \xymatrix{
      \epsilon \circ H^0(\bM) \circ H^0(\bM) \ar@{=}[r]
      \ar[d]_{\epsilon(H^0(\dmu))} & M \circ M \circ \epsilon
      \ar[d]^{\mu(\epsilon)} \ar@{}[drr]|{\text{and}} & & \epsilon \circ
      H^0(\idf{\dC}) \ar@{=}[r] \ar[d]_{\epsilon(H^0(\deta))} & \idf{\cC}
      \circ \epsilon \ar[d]^{\eta(\epsilon)} \\
      \epsilon \circ H^0(\bM) \ar@{=}[r] & M \circ \epsilon, & & \epsilon
      \circ H^0(\bM) \ar@{=}[r] & M \circ \epsilon.
    }
  \end{equation*}
  \label{dfn:liftdgmd}
\end{definition}

The definition of DG monads which we get from the theory of enriched
categories is not sufficient for our purposes. We use the following
definition of modules over DG monads.

\begin{definition}
  Let $\dC$ be a DG category and $(\bM, \dmu, \deta)$ be a DG monad on
  $\dC$. A DG $\bM$-module is a pair $(C, \dlambda)$ consisting of an object
  $C$ in $\dC$ and a morphism $\dlambda : \bM C \to C$ such that
  \begin{enumerate}
    \item $\deg \dlambda = 0$,
    \item $d \dlambda = 0$, and
    \item the following diagrams commute.
      \begin{equation*}
        \xymatrix{
          \bM^2 C \ar[r]^{M \dlambda} \ar[d]_{\dmu_C} & \bM C
          \ar[d]^{\dlambda} & & & C \ar[r]^{\deta_C} \ar@{=}[dr] & \bM C
          \ar[d]^{\dlambda} \\
          \bM C \ar[r]^{\dlambda} & C & & & & C
        }
      \end{equation*}
  \end{enumerate}
  As in the classical case, a morphism $\varphi : (C, \dlambda_C) \to (D,
  \dlambda_D)$ is a morphism $\varphi : C \to D$ in $\dC$ such that the
  following diagram commutes.
  \begin{equation}
  \xymatrix{
    \bM C \ar[r]^{\bM \varphi} \ar[d]_{\dlambda_C} &
    \bM D \ar[d]^{\dlambda_D} \\
    C \ar[r]^{\varphi} &
    D
  }
  \label{dia:morpdmod}
  \end{equation}
  The category of DG $\bM$-modules will be denoted by $\moduledg{\bM}$.
\end{definition}


\begin{remark}
  Note that in the above definition, 3 implies 1, since definition
  \ref{dfn:dgnattrn} demands $\deta$ and $\dmu$ to be closed of degree
  $0$. Thus the only extra condition here is closedness of $\dlambda$. The
  closedness of $\dlambda$ implies that for any $\varphi : C \to D$ fitting
  into a diagram like \eqref{dia:morpdmod}, $d \varphi$ is also a morphism in
  $\moduledg{\bM}$, making $\moduledg{\bM}$ a DG category.

  Now it is easy to check that $\bF_{\bM} : \dC \to \moduledg{\bM}$ defined
  by
  \begin{align*}
  \bF_{\bM}(C) &= (\bM C, \dmu_C) &
  \bF_{\bM}(\varphi) &= \bM \varphi, \\
  \bG_{\bM}((C, \dlambda_C)) &= C &
  \bG_{\bM}(\varphi) &= \varphi.
  \end{align*}
  defines a DG adjunction
  \begin{equation*}
  \adjunct{\bF_{\bM}}{\dC}{\moduledg{\bM}}{\bG}.
  \end{equation*}

  Also it is easy to prove a lemma similar to lemma \ref{lem:emenrmnd}.
  Given any adjoint pair $\adjunct{\bF}{\dC}{\dD}{\bG}$, let $\bM = \bG
  \circ \bF$. Following the lines of the proof of \ref{lem:emenrmnd}, one
  can define $\bL : \dD \to \moduledg{\bM}$, by $\bL D = (\bG D, \dlambda)$
  where $\dlambda : \bM \bG D \to \bG D$ is the map
  \begin{equation*}
  \bM \bG D =  (\bG \bF) \bG D =  \bG (\bF \bG ) \bD \xrightarrow{\bG
  \depsilon_D} D,
  \end{equation*}
  $\depsilon_D$ being the counit natural transformation $\bF \bG \to
  \idf{\dD}$. It is obvious what $\bL$ should be on morphisms.
  \label{rmk:dgcatmmd}
\end{remark}

\begin{lemma}
  If $\dC$ is a strongly pretriangulated DG category and $(\bM, \dmu,
  \deta)$ is a DG monad on $\dC$, then $\moduledg{\bM}$ is a strongly
  pretriangulated DG category. Thus, $H^0(\moduledg{\bM})$ is a triangulated
  category.
  \label{lem:dgmmodtr}
\end{lemma}

\begin{proof}
  By remark \ref{rmk:dgcatmmd}, $\moduledg{\bM}$ is a DG category.  To show
  that it is strongly pretriangulated, we have to show that every object has
  a shift and every closed, degree $0$ morphism has a cone (see remark
  \ref{rmk:altdscpt}). We define these as follows.
  \begin{itemize}
    \item For $\underline{A} = (A, \dlambda_A)$, set $\underline{A}[1] =
      (A[1], \dlambda[1])$.
    \item For $\underline{A} = (A, \dlambda_A)$ and $\underline{B} = (B,
      \dlambda_B)$, and $\varphi : \underline{A} \to \underline{B}$ a closed
      degree zero morphism, then we define $\Cone(\varphi)$ as follows. By
      abuse of notation let $\varphi : A \to B$ in $\dC$ be the underlying
      morphism. Let $C = \Cone(\varphi)$ in $\dC$. Thus we have a diagram of
      the form
      \begin{equation*}
        \xymatrix{
          \bM A \ar[r]^{\bM \varphi} \ar[d]_{\lambda_A} &
          \bM B \ar[r] \ar[d]^{\lambda_B} &
          \bM C \ar@{.>}[d]^{\lambda_C} \ar@{=}[r] &
          \Cone(\bM \varphi)
          \\
          A \ar[r] &
          B \ar[r] &
          C &
        }
      \end{equation*}
    By construction of cone, there exists a morphism $\lambda_C : \bM C
    \to C$ making the above diagram commutative (see \ref{prp:funccone}).
    Using functoriality of cones and $\bM$, it is easy to check that
    $\underline{C} = (C, \lambda_C)$ is DG isomorphic to the cone of
    $\varphi$ in $\moduledg{\bM}$.
  \end{itemize}
  The last assertion of the lemma follows from proposition
  \ref{prp:h0pretrtr}.
\end{proof}

\begin{definition}
  Suppose we have a triangulated category $\cC$, a monad $(M, \mu, \eta)$ on
  $\cC$, a strong enhancement $(\dC, \epsilon)$ of $\cC$ and a lift $\bM$ of
  $M$ to $\dC$. We shall denote the triangulated category
  $H^0(\moduledg{\bM})$ by $\Dtr{\dC, \bM}$, $\Dtr{\bM}$ or just $\Dtr{}$.
  \label{dfn:dtrMmods}
\end{definition}

\begin{proposition}
  Suppose $(M, \mu, \eta)$ be an exact monad on a triangulated category
  $\cC$. Assume that there exists an enhancement $(\dC, \epsilon)$ of $\cC$
  and a lift $\bM$ of $M$ to $\dC$. Then there exists a pair of exact
  adjoint functors
  \begin{equation*}
    \adjunct{F^{tr}}{\cC}{\Dtr{\bM}}{G^{tr}}
  \end{equation*}
  satisfying $M = G^{tr} \circ F^{tr}$.
  \label{prp:factexad}
\end{proposition}

\begin{remark}
  Note that this partially answers Balmers question in \cite[remark
  2.9]{balmer:septrcat}, where he asks if given a monad $M$ on a
  triangulated category $\cC$, there exists a triangulated category
  $\cD$ and an exact adjunction $\adjunct{F}{\cC}{\cD}{G}$ such that
  $G \circ F = M$.
\end{remark}

\begin{proof}[Proof of proposition \ref{prp:factexad}]
  First assume that the lift $\dC$ is strongly pretriangulated. Then by
  lemma \ref{lem:dgmmodtr}, $\moduledg{\bM}$ is also strongly
  pretriangulated. As mentioned in remark \ref{rmk:dgcatmmd}, given $\dC$
  and $\bM$, one can construct an adjunction
  \begin{equation*}
  \adjunct{\bF_{\bM}}{\dC}{\moduledg{\bM}}{\bG_{\bM}}
  \end{equation*}
  Now the required adjunction mentioned above, is obtained by applying the
  2-functor $H^0$ (see \ref{lem:h0strict}) to the above adjunction and
  setting $F^{tr} = H^0(\bF_{\bM})$, and $G^{tr} = H^0(\bG_{\bM})$. Note,
  that $\cC = H^0(\dC)$ follows from the definition of enhancement
  (definition \ref{dfn:liftdgmd}), and $\Dtr{\bM} = H^0(\moduledg{\bM})$ is
  by definition \ref{dfn:dtrMmods}.

  Now for general lifts $\dC$, let $\boldi : \dC \to \pretr{\dC}$ be the
  canonical inclusion functor. By lemma \ref{lem:pretr2fn},
  $(\pretr{\bM}, \pretr{\dmu}, \pretr{\deta})$ is a DG monad on
  $\pretr{\dC}$. It is easy to see that $\pretr{\bM}$ is also a lift of
  $M$, as $H^0(\boldi)$ is an equivalence of categories. Thus, the same
  argument as in the first paragraph applied on $\pretr{\dC}$ and
  $\pretr{\bM}$ in place of $\dC$ and $\bM$ respectively gives us an
  adjunction
  \begin{equation*}
  \adjunct{\bF_{\pretr{\bM}}}{\pretr{\dC}}{\moduledg{\pretr{\bM}}}{\bG_{\pretr{\bM}}}
  \end{equation*}
  Now taking $H^0$ as before we get the required result.
\end{proof}

\tododone{The following seems to be Theorem 7.10(3). So I am removing it.
}

\begin{remark}
  In proposition \ref{prp:factexad}, if $G^{tr}$ is separable, then
  \cite[Main Theorem 5.17(d)]{balmer:septrcat} implies that $\Dtr{\bM}$ is
  equivalent to $\module{M}$ and in this case $\module{M}$ is triangulated.
\end{remark}

\tododone{Include the fact that if $G^{tr}$ is separable, Balmer's proposition
5.17 in the paper on separable triangulated category implies that
$\Dtr{\bM}$ is equivalent to $\module{M}$ and hence $\module{M}$ is
triangulated.}


\tododone{Insert example 3.7 and lemma 3.8 from the old version.}

We conclude this subsection by studying a couple of examples.

\begin{example}
  Consider an abelian category $\cA$. Let $\cC^{\#}(\cA)$ be the additive
  category of complexes in $\cA$. Here $\#$ can be one of $+$, $-$, $b$ or
  it can be empty depending on whether we look at bounded above, bounded
  below, bounded or unbounded complexes. Let $\cK^{\#}(\cA)$ be the
  corresponding homotopy category. $\cC_{dg}^{\#}(\cA)$ will denote the DG
  category of complexes in $\cA$ (see \cite{keller:derdgcat}, for example).
  Suppose $(M, \mu, \eta)$ is an exact monad on $\cA$.

  With this setup, let $DG : \abcat \to \dgcat$ be the functor from the
  category of small abelian categories to small DG categories defined by
  $DG(\cA) = \cC_{dg}^{\#}(\cA)$ on objects of $\abcat$, $DG(F : \cA \to
  \cB) = F^{\bullet} : \cC_{dg}^{\#}(\cA) \to \cC_{dg}^{\#}(\cA)$ given by
  $F^{\bullet} (A^{\bullet}, d^{\bullet}) = (F(A^{\bullet}),
  F(d^{\bullet}))$. For natural transformations $\alpha: F \to G$ between
  functors in $\abcat$, one can define a natural transformation $DG(\alpha)
  : DG(F) \to DG(G)$ as $DG(\alpha)_{(A^{\bullet}, d^{\bullet})} :
  DG(F)(A^{\bullet}, d^{\bullet}) \to DG(G)(A^{\bullet}, d^{\bullet})$ as
  the collection $\alpha^{\bullet} : F(A^{\bullet}) \to G(A^{\bullet})$.
  With this definition it is easy to see that $DG$ is a 2-functor. Let us
  denote the monad $(DG(M), DG(\mu), DG(\eta))$ by $(\bM_{dg},
  \dmu_{dg}, \deta_{dg})$. This is a natural example of a DG monad.
  \label{exl:monadabc}
\end{example}

\begin{lemma}
  $\moduledg{\bM_{dg}}$ is equivalent to $\cC^{\#}_{dg}(\module{M})$. Thus,
  $H^0(\moduledg{\bM_{dg}})$ is equivalent to $\cK^{\#}(\module{M})$.
  \label{lem:egcomplx}
\end{lemma}

\begin{proof}
  The elements of $\cC^{\#}_{dg}(\module{M})$ are complexes $(x^{\bullet},
  \lambda^{\bullet}), d^{\bullet}$ of objects in $\module{M}$. By
  definition, of morphism of modules, $(x^{\bullet}, d^{\bullet})$ is a
  complex and an element of $\cC^{\#}(\cA)$. It is easy to check that
  $((x^{\bullet}, d^{\bullet}), \lambda^{\bullet})$ is an element of
  $\moduledg{\bM_{dg}}$. The functor $\Psi : \cC^{\#}_{dg}(\module{M}) \to
  \moduledg{\bM_{dg}}$ is easily shown to be an equivalence.

  The second part of the sentence is obtained by applying $H^0$ to the above
  equivalence.
\end{proof}

\subsection{Weak monads on DG categories} \label{ssc:weakmond}
\tododone{2.21 -- 2.32}

\tododone{weak lifts, 2.21 - 2.32, natural transformations lift to weak natural
transformations.}

\tododone{relate to lemma 3.7, when does a monad on K come from Cdg? Any monad
lifts to weak monad. ???}

In this subsection, we shall weaken the definition of a natural
transformation between two DG functors.

\begin{definition}
  Suppose $\bF$ and $\bG$ are two DG functors between DG categories $\dC$
  and $\dD$. A \emph{weak natural transformation} $\dalpha : \bF \to \bG$ is
  a collection of morphisms $\set{\dalpha_c : \bF c \to \bG c}{c \text{ is an
    object in } \dC}$, such that
  \begin{itemize}
    \item for all objects $c$ in $\dC$, $\dalpha_c$ has degree zero;
    \item satisfies $d \dalpha_c = 0$; and
    \item the induced maps $H^0(\dalpha_c) \in H^0(\Hom_{\dD}(\bF c, \bG c))$
      gives a natural transformation $H^0(\dalpha) : H^0(\bF) \to H^0(\bG)$.
  \end{itemize}

  A weak natural transformation $\dalpha$ is said to be a \emph{weak natural
  isomorphism} if $H^0(\dalpha)$ is a natural isomorphism.

  A \emph{weak inverse}, $\dbeta$ of a weak natural isomorphism $\dalpha$
  is a choice of morphisms $\dbeta_x : \bG x \to \bF x$, for each object
  $x$ in $\dC$, such that 
  \begin{itemize}
    \item $\deg \dbeta_x = 0$,
    \item $d \dbeta_x = 0$
  \end{itemize}
  for all objects $x$ in $\dC$, and $x \mapsto H^0(\dbeta_x)$ is a natural
  transformation from $H^0(\bG)$ to $H^0(\bF)$, and $H^0(\dbeta) =
  H^0(\dalpha)^{-1}$ as natural transformations on additive categories.
  \label{def:everweak}
\end{definition}

\begin{definition}
  Suppose $\dC$ is a DG category and $\bM : \dC \to \dC$ is a DG
  endofunctor. Let $\dmu : \bM^2 \to \bM$ and $\deta : \id \to \bM$ be weak
  DG natural transformations. We say that $(\bM, \dmu, \deta)$ is a
  \emph{weak DG monad} on $\dC$ if $(H^0(\bM), H^0(\dmu), H^0(\deta))$ is a
  monad on the additive category $H^0(\dC)$.
  \label{dfn:wkmonads}
\end{definition}

Now we define weak modules over weak monads.

\begin{definition}
  Let $\dC$ be a DG category and $(\bM, \dmu, \deta)$ be a weak monad on
  $\dC$. An \emph{weak $\bM$-module} is a pair $(x, \dlambda)$ where
  \begin{itemize}
    \item $x$ is an object in $\dC$,
    \item $\dlambda : \bM x \to x$ is a morphism of degree $0$ such that
      $d \dlambda = 0$, and
    \item $\dlambda \circ \bM \dlambda - \dlambda \circ \dmu_x = d \drho$
      where $\drho : \bM^2 x \to x$ is a morphism of degree $-1$.
    \item $\dlambda \circ \deta_x - \idf{x} = d \dnu$ where $\dnu : x \to x$
      is a morphism of degree $-1$.
  \end{itemize}
\end{definition}

\begin{remark}
  In the above definition, $d(\dlambda \circ \bM \dlambda - \dlambda \circ
  \dmu_x)$ can easily be checked to be $0$ using Liebnitz rule for $d$.
  Similarly, $d(\dlambda \circ \deta_x - \idf{x}) = 0$. Thus the last two
  conditions can be interpreted as $H^0(\dlambda \circ \bM \dlambda -
  \dlambda \circ \dmu_x) = H^0(\dlambda \circ \deta_x - \idf{x}) = 0$. 
\end{remark}

\begin{definition}
  Suppose $\dC$ is a DG category and $(\bM, \dmu, \deta)$ be a weak monad on
  it. Suppose $(x, \dlambda)$ and $(x', \dlambda')$ are two weak
  $\bM$-modules. A \emph{weak $\bM$-module morphism} $\varphi : (x,
  \dlambda) \to (x', \dlambda')$ is a morphism $\varphi : x \to x'$ in
  $\Hom_{\dC}(x, x')$ such that
  \begin{equation*}
    \dlambda' \circ \bM \varphi - \varphi \circ \dlambda =
    \sum_{i=1}^{n} f_i \circ d \theta_i
  \end{equation*}
  for some finite set of morphisms $f_i$ and $\theta_i$ in $\dC$.
  \label{dfn:wkmodmrp}
\end{definition}

\begin{remark}
  Asking for $\dlambda' \circ \bM \varphi - \varphi \circ \dlambda = d
  \theta$ in the above definition, is not enough for the resulting category
  of modules to be closed under composition.
\end{remark}

\begin{proposition}
  The weak $\bM$-modules along with the weak $\bM$-module morphisms form a
  DG category, which we denote by $\moduleweak{\bM}$. There exists an
  essentially surjective functor $H^0(\moduleweak{\bM}) \to \module{M}$
  where $M = H^0(\bM)$.
\end{proposition}

\begin{proof}
  That weak $\bM$ modules form a DG category is a routine computation.
  Consider the functor
  \begin{equation*}
    \Phi : H^0(\moduleweak{\bM}) \to \module{M}
  \end{equation*}
  defined by $\Phi((x, \dlambda)) = (x, \lambda)$, where $\lambda =
  H^0(\dlambda)$ is the image of $\dlambda$ in $\Hom_{H^0(\dC)}(Mx, x) =
  H^0(\Hom_{\dC}(\bM x, x)$. Suppose $f : (x, \dlambda) \to (x', \dlambda')$
  is a morphism in $H^0(\moduleweak{\bM})$. Let $\mathbf{f} : (x, \dlambda)
  \to (x', \dlambda')$ be any closed morphism in $\moduleweak{\bM}$ whose
  $H^0$ is $f$. Then $\Phi(f) = H^0(\mathbf{f})$ in $H^0(\dC)$ satisfies
  $\lambda' \circ M \Phi(f) = \Phi(f) \circ \lambda$, and thus defines a
  morphism in $\module{M}$. Since $\Hom_{\moduleweak{\bM}}((x, \dlambda),
  (x', \dlambda')) \subset \Hom_{\dC}(x, x')$, $\Phi(f)$ is well defined.

  That $\Phi$ is essentially surjective, is obvious by definition.
\end{proof}

Finally, to make sense of the Eilenberg-Moore construction in this setting
we need a notion of weak adjoints.

\begin{definition}
  Two DG functors $\bF : \dC \to \dD$ and $\bG : \dD \to \dC$ between DG
  categories $\dC$ and $\dD$ are said to be \emph{weak adjoint} if
  $\adjunct{H^0(\bF)}{H^0(\dC)}{H^0(\dD)}{H^0(\bG)}$ is an adjunction
  between additive categories.
\end{definition}

\begin{lemma}
  If $\adjunct{\bF}{\dC}{\dD}{\bG}$ is a weak adjoint pair of functors between
  DG categories, then $\bG \circ \bF$ defines a weak DG monad on $\dC$.
\end{lemma}
\tododone{write $\dmu$ and $\deta$}
\begin{proof}
  Recall that in the additive category case the multiplication and unit
  natural tranformations of a monad are induced by the unit and the co-unit
  natural transformations of the adjoint functors. We do the same thing
  here. We choose any weak natural transformation $\depsilon : \bF \circ \bG
  \to \idf{\dC}$ by choosing $\depsilon_A : \bF \circ \bG(A) \to A$ to be a
  degree zero closed morphism such that $H^0(\depsilon_A) = \varepsilon_A :
  H^0(\bF) \circ H^0(\bG) (A) \to A$ where $\varepsilon_A$ is the co-unit
  for the adjunction $\adjunct{H^0(\bF)}{H^0(\dC)}{H^0(\dD)}{H^0(\bG)}$.

  Thus, we can define $\dmu_A$ to be $\bG \depsilon_A \bF$. From the
  construction, it is clear that $H^0(\dmu)$ is a natural transformation.
  One chooses $\deta_A : A \to \bM A = \bG \circ \bF A$ to be any closed
  morphism such that $H^0(\deta_A) : A \to H^0(\bG) \circ H^0(\bF)$ is the
  unit natural transformation for the adjoint pair $(\bF, \bG)$. Thus
  $\deta_A$ is a weak natural transformation.

  It is obvious that $(\bM, \dmu, \deta)$ is a weak monad.
\end{proof}

\begin{lemma}
  Let $\dC$ and $(\bM, \dmu, \deta)$ be a DG category and a weak monad on it
  as above. Then one has a weak adjunction
  \begin{equation*}
    \adjunct{\bF_{\bM}}{\dC}{\moduleweak{\bM}}{\bG_{\bM}},
  \end{equation*}
  where $\bG_{\bM}$ is the forgetful functor $\bG_{\bM}((x, \lambda)) = x$
  and $\bF_{\bM}(x) = (\bM x, \dmu_x)$. Furthermore, $\bM = \bG_{\bM} \circ
  \bF_{\bM}$.
\end{lemma}

\begin{proof}
  It is an easy exercise to check that $\bF_{\bM}$ and $\bG_{\bM}$ are
  functors (see \ref{dfn:wkmodmrp}). The proof of the fact that $\bF_{\bM}$
  is a left adjoint to $\bG_{\bM}$ is similar to the classical case; see for
  example \cite[chapter VI, section 2, theorem 1]{maclane:catworkmath} and
  \cite{eilenbergmoore:adjfunctrip}.
\end{proof}

\begin{definition}
  The category of weak free $\bM$-modules is the full subcategory of
  $\moduleweak{\bM}$ generated by the image of $F_{\bM}$. It will be denoted
  by $\freeweak{\bM}$.
\end{definition}

To make sense of the universal properties of $\moduleweak{\bM}$ and
$\freeweak{\bM}$, we need the following definition.

\begin{definition}
  A DG functor $\bH : \dA \to \dB$ between two DG categories, satisfying
  certain properties, is said to be \emph{unique up to a quasi-unique weak
  natural isomorphism}, if for any other DG functor $\bH' : \dA \to \dB$
  satisfying the same properties, there exists a weak natural isomorphism
  $\dalpha : \bH \to \bH'$ such that $H^0(\dalpha)$ has a unique choice.
\end{definition}

\begin{proposition}
  Let $\dC$ and $\dD$ be DG categories and $\adjunct{\bF}{\dC}{\dD}{\bG}$ be a
  weak adjoint pair of functors. Let $\bM = \bG \circ \bF$ be the corresponding
  weak DG monad on $\dC$. Then we have two DG functors:
  \begin{itemize}
    \item $\bK : \dD \to \moduleweak{\bM}$ and
    \item $\bL : \freeweak{\bM} \to \dD$,
  \end{itemize}
  both unique up to a quasi-unique weak natural isomorphism, such that the
  arrows fit into the following diagram:
  \begin{equation*}
    \xymatrix{
    & & \dC \ar@<1pt>@^{->}[ddll]^{\bF_{\bM}}
    \ar@<1pt>@^{->}[ddrr]^{\bF_{\bM}} \ar@<1pt>@^{->}[dd]^{\bF} & & \\
      & & & & \\
    \freeweak{\bM} \ar@<1pt>@^{->}[uurr]^{\bG_{\bM}} \ar[rr]_-{\bL} & &
    \dD \ar@<1pt>@^{->}[uu]^{\bG} \ar[rr]_-{\bK} & & \moduleweak{\bM}
    \ar@<1pt>@^{->}[uull]^{\bG_{\bM}} 
    }
  \end{equation*}
  where $\bL \circ \bF_{\bM} = \bF$, $\bG \circ \bL = \bG_{\bM}$, $\bF_{\bM}
  = \bK \circ \bF$ and $\bG = \bG_{\bM} \circ \bK$.
\end{proposition}

\begin{proof}
  If $\depsilon : \bF \circ \bG \to \idf{\dC}$ is the co-unit weak natural
  transformation, one defines $\bL : \moduleweak{\bM} \to \dD$ by $\bL((\bM
  x, \dmu_x)) = \depsilon_{\bF x} (\bF (\bM x))$. This is easily checked to
  be well defined in lines of the proof of \cite[chapter VI, section 2,
  theorem 1]{maclane:catworkmath}. One can also define $\bK D = (\bG D, \bG
  \depsilon_D)$. Now the theorem follows from standard arguments about
  adjoint pairs of functors (see loc.~cit.).
\end{proof}

\begin{theorem}
  Let $(M, \mu, \eta)$ be a monad on a triangulated category $\cC$.
  Suppose $\cC$ admits an enhancement $(\dC, \epsilon)$, and an endofunctor
  $\bM : \dC \to \dC$ which satisfies $H^0(\bM) = M$. Then there exists a
  triangulated category $\cD$ and an exact adjunction
  $\adjunct{F}{\cC}{\cD}{G}$ such that $G \circ F = M$.
\end{theorem}

\begin{proof}
  The proof is in the same lines as that of proposition \ref{prp:factexad}.
  Choose degree 0 closed natural transformations $\dmu : \bM^2 \to \bM$ and
  $\deta : \idf{\dC} \to \bM$ such that $H^0(\dmu) = \mu$ and $H^0(\deta) =
  \eta$. By definition, $(\bM, \dmu, \deta)$ is a weak monad on $\dC$.
  Consider the weak adjunction
  \begin{equation*}
  \adjunct{\bF_{\bM}}{\dC}{\moduleweak{\bM}}{\bG_{\bM}}.
  \end{equation*}
  Applying the strict 2-functor $\PreTr$ we get a weak adjunction
  \begin{equation*}
  \adjunct{\pretr{\bF_{\bM}}}{\pretr{\dC}}{\pretr{\moduleweak{\bM}}}{\pretr{\bG_{\bM}}}.
  \end{equation*}
  Now apply the strict 2-functor $H^0$ to get an adjunction
  \begin{equation*}
  \adjunct{F}{\cC}{\dD}{G}
  \end{equation*}
  where,
  \begin{align*}
  F &= H^0(\pretr{\bF_{\bM}}) &
  G &= H^0(\pretr{\bG_{\bM}}) \\
  \cC &= H^0(\dC) \cong H^0(\pretr{\dC}) &
  \cD &= H^0(\pretr{\moduleweak{\bM}}).
  \end{align*}
  Clearly, $\cD$ is triangulated and $F$ and $G$ are exact.
\end{proof}

\section{Bousfield localizations in DG categories} \label{sec:bousfield}
\tododone{3.2 -- 3.5}

\begin{definition}
  We shall call a DG endofunctor $\bL : \dC \to \dC$ on a DG category
  $\dC$, a \emph{weak Bousfield localization} functor if there is a weak
  natural transformation $\deta : \idf{\dC} \to \bL$ such that $\bL \deta :
  \bL \to \bL^2$ is a weak natural isomorphism and $\bL \deta = \deta \bL$.
  \label{def:boulocdg}
\end{definition}

\begin{definition}
  By a Bousfield localization functor on a triangulated category $\cT$, we
  mean an exact functor $L : \cT \to \cT$ along with a natural
  transformation $\eta : \id_{\cT} \to L$ such that $L \eta : L \to L^2$ is
  a natural isomorphism and $L \eta = \eta L$. For details, refer to
  \cite{krause:localization}.
\end{definition}

\begin{lemma}
  Suppose $\dC$ is a pre-triangulated DG category. $\bL : \dC \to \dC$ being
  a weak Bousfield localization functor is equivalent to  $L = H^0(\bL)$
  being a Bousfield localization functor on the triangulated category $\cC =
  H^0(\dC)$.
  \label{lem:h0wkblbl}
\end{lemma}

\begin{proof}
  This is by definition of weak natural transformations.
\end{proof}

\begin{definition}
  We call a weak monad $(\bM, \dmu, \deta)$ on a pre-triangulated DG
  category $\dC$ to be \emph{separable} if $(H^0(\bM), H^0(\dmu),
  H^0(\deta))$ is separable (see \ref{dfn:sepmonad}) as a monad on the
  triangulated category $H^0(\dC)$.
\end{definition}

\begin{lemma}
  Let $\dC$ be a DG category and $\bL : \dC \to \dC$ be a weak Bousfield
  localization functor. Let $\deta : \idf{\dC} \to \bL$ be the weak natural
  transformation as in definition \ref{def:boulocdg}. Let $\dmu : \bL^2 \to
  \bL$ be a weak inverse (see \ref{def:everweak}) of $\bL \deta$. Then
  $(\bL, \dmu, \deta)$ is a separable weak monad.
\end{lemma}

\begin{proof}
  This is obvious as $L = H^0(\bL)$ is a Bousfield localization functor on
  $\cC = H^0(\dC)$. If $\eta = H^0(\deta)$, it is easy to see that
  $(L, \mu = (L \eta)^{-1}, \eta)$ is a monad on $\cC$ with $\mu$
  invertible. Set $\sigma = L \eta$. To prove that $L$ is separable, we need
  \begin{equation*}
  L (L \eta)^{-1} \circ (L \eta) L = (L \eta)^{-1} L \circ L (L \eta) = L
  \eta \circ (L \eta)^{-1} = \idf{L^2}.
  \end{equation*}
  which is easy to see as $L \eta = \eta L$.
\end{proof}

\begin{definition}
  A full DG subcategory $\dE$ of a DG category $\dC$ is said to be
  \emph{weakly admissible} if the inclusion DG functor $\dE \to \dC$ has a
  weak right adjoint.
\end{definition}

We define the kernel of a functor, and prove a criteria for a functor to be
a weak Bousfield localization functor. This is a DG version of
\cite[proposition 4.9.1]{krause:localization}.

\begin{definition}
  Let $\dC$ be a pre-triangulated DG category and $\bL : \dC \to \dC$ be a
  pre-exact DG endofunctor. The \emph{kernel of $\bL$}, $\ker \bL$ is the
  full subcategory of $\dC$ consisting of objects $c$ such that $\idf{\bL c}
  = d \theta$ for some morphism $\theta : \bL c \to \bL c$ of degree $-1$.
\end{definition}

\begin{proposition}
  Let $\dC$ be a pre-triangulated DG category and $\bL : \dC \to \dC$ be a
  pre-exact DG endofunctor. Then the following are equivalent:
  \begin{enumerate}
    \item $\bL$ is a weak Bousfield localization functor.
    \item $\ker \bL$ is a weakly admissible subcategory of $\dC$.
    \item The Drinfeld quotient $\dC \to \dC / (\ker \bL)$ admits a weak
      right adjoint.
  \end{enumerate}
  \label{prp:bouslocw}
\end{proposition}

\begin{proof}
  Note that $\cC = H^0(\dC)$ is triangulated and by lemma
  \ref{lem:h0wkblbl}, $L = H^0(\bL)$ is a Bousfield localization functor on
  $\cC$. 
  
  Moreover, $c$ belongs to $\ker \bL$, if and only if $H^0(\idf{\bL c}) =
  \idf{Lc} = 0 \in \Hom_{\cC}(Lc, Lc) = H^0(\Hom_{\dC}(\bL c, \bL c)) =
  \{0\}$; that is, $Lc \cong 0$ in $\cC$. Thus $c$ is an object of $\ker L$
  in $\cC$. Note that the arguments are reversible and hence $c$ belongs to
  $\ker \bL$ in $\dC$ if and only if $c$ belongs to $\ker L$ in $\cC$.
  
  Furthermore, $\ker \bL$ is a weakly admissible subcategory when and only
  when $\ker L$ is an admissible subcategory of $\cC$.
  
  Finally, the exitence of weak right adjoint of the Drinfeld quotient $\dC
  \to \dC / (\ker \bL)$ is equivalent to the exitence of a right adjoint to
  the Verdier localization functor $\cC \to \cC / (\ker L)$.

  Now the theorem is a direct consequence of \cite[proposition
  4.9.1]{krause:localization}.
\end{proof}

\tododone{Remark 3.6 fails as we really do not have H0(M-hmod) = M-mod. We need
a different statement here.}

\tododone{Write the observation that $\module{L}$ is equivalent to $\cC / \ker
L$ and has an enhancement $(\ker {\bL})^{\perp}$.}

\begin{example}
  Suppose $\bL : \dC \to \dC$ is a weak Bousfield localization functor on a
  pretriangulated DG category, as in proposition \ref{prp:bouslocw}. Suppose
  $\dK = \ker \bL$. Let $\cC = H^0(\dC)$, $L = H^0(\bL)$ and $\cK = \ker L$.
  It is easy to see that $H^0(\dK) = \cK$.  Thus $L : \cC \to \cC$ is a
  Bousfield localization functor. If $Q : \cC \to \cC / \cK$ is the Verdier
  quotient functor, then from the proof of theorem 4.9.1 in
  \cite{krause:localization}, we have the following facts.
  \begin{enumerate}
  \item $Q$ has a right adjoint $Q_{\rho} : \cC / \cK \to \cC$.
  \item $Q_{\rho} \circ Q = L$ and $Q \circ Q_{\rho}
    \xrightarrow[\cong]{\epsilon} \id_{\cC / \cK}$.
  \item $Q_{\rho}$ is fully faithful.
  \end{enumerate}
  Thus the co-unit natural transformation $\epsilon$ is invertible, and
  hence has a section. By Balmer \cite[remark 3.9]{balmer:septrcat} (see
  also Rafael \cite[theorem 1.2]{rafael:sepfunct}), $Q_{\rho}$ is separable,
  and hence stably separable. Thus, one can apply theorem
  \ref{thm:Mmodntrg}(d) to conclude that in this case $\module{L}$ is
  triangulated and 
  \begin{equation*}
  \cK^{\perp} \cong \cC / \cK \cong \module{L},
  \end{equation*}
  where $\cong$ is used to denote an equivalence of triangulated categories.
  (cf.~\cite[example 6.3]{balmer:septrcat}.)

  From the above discussion it is also clear that $\module{L}$ admits an
  enhancement $\dK^{\perp} = (\ker \bL)^{\perp}$ where the equivalence
  $H^0(\varepsilon) : H^0(\dK^{\perp}) \to \cK^{\perp}$ is induced by the
  equivalence between $H^0(\dC)$ and $\cC$.
\end{example}

\section{$G$-equivariant derived categories} \label{sec:gequivdc}
\tododone{2.16 -- 2.20, drop 2.15, 3.1 : elaborate on the last part, may be
break it off as a proposition}

\tododone{Refer to the definition of $\perf{\dC}$ in section 3.2}

\tododone{Break off the following into three definiitons: (Mod), Q and Gamma}

\subsection{A general construction} \label{ssc:gequivgen}

\tododone{Clarify that all the categories are 2-categories and all the functors
are 2-functors}

\begin{definition}
  Let $\cA'$ be the following category.

  The objects of $\cA'$ are diagrams of the form
  \begin{equation*}
    \xymatrix{
      \dC \ar[r] & \perf{\dC} \ar@^{->}@<1pt>[d]^{\bF} \\
      & \dD \ar@^{->}@<1pt>[u]^{\bG}
    }
  \end{equation*}
  where $\dC$ and $\dD$ are DG categories, $\perf{\dC}$ denotes the DG
  category of perfect DG $\dC$-modules (see definition \ref{dfn:perfdgct})
  and $\adjunct{\bF}{\perf{\dC}}{\dD}{\bG}$ is an adjunction.

  Given two objects
  \begin{equation*}
    \xymatrix{
      \dC \ar[r] & \perf{\dC} \ar@^{->}@<1pt>[d]^{\bF} &
      \ar@{}[dr]|{\text{\normalsize{and}}} & &
      \dC' \ar[r] & \perf{(\dC')} \ar@^{->}@<1pt>[d]^{\bF'} \\
      & \dD \ar@^{->}@<1pt>[u]^{\bG} & & &
      & \dD' \ar@^{->}@<1pt>[u]^{\bG'}
    }
  \end{equation*}
  in $\cA'$, a morphism $(\bC, \bD)$ is a pair of functors
  $\bC : \dC \to \dC'$ and $\bD : \dD \to \dD'$ such that the following
  diagrams of functors commute
  \begin{equation*}
    \xymatrix{
      \perf{\dC} \ar[r]^{\perf{\bC}} \ar[d]_{\bF} & \perf{(\dC')}
      \ar[d]^{\bF'} & & & \perf{\dC} \ar[r]^{\perf{\bC}} & \perf{(\dC')} \\
      \dD \ar[r]^{\bD} & \dD' & & & \dD \ar[r]^{\bD} \ar[u]^{\bG} & \dD'
      \ar[u]_{\bG'},
    }
  \end{equation*}
  or in other words the relevant squares in the following diagram commutes
  \begin{equation*}
    \xymatrix{
      & \dC \ar[r] \ar[dl]_{\bC} & \perf{\dC} \ar@^{->}@<1pt>[d]^{\bF}
      \ar[dl]_{\perf{\bC}} \\
      \dC' \ar[r] & \perf{(\dC')} \ar@^{->}@<1pt>[d]^{\bF'} & \dD
      \ar@^{->}@<1pt>[u]^{\bG} \ar[dl]^{\bD} \\
      & \dD'. \ar@^{->}@<1pt>[u]^{\bG'} &
    }
  \end{equation*}
  Composition is given by $(\bC', \bD') \circ (\bC, \bD) = (\bC' \circ \bC,
  \bD' \circ \bD)$.  It is easy to see that this forms a category.
\end{definition}

\begin{definition}
  For an object $\objAline{\dC}{\dD}{\bF}{\bG}$ in $\cA'$, we define 
  \begin{equation*}
    \dQ \big( \objAline{\dC}{\dD}{\bF}{\bG} \big)
  \end{equation*}
  to be the full subcategory of $\dD$ consisting of all objects $d$ such
  that $\bG d$ is homotopically equivalent (see \ref{def:homequiv}) to the
  image of an object in $\dC$.

  Let
  \begin{equation*}
    A = \objA{\dC}{\dD}{\bF}{\bG} \quad \text{and} \quad A' =
    \objAb{\dC'}{\dD'}{\bF'}{\bG'}
  \end{equation*}
  be two objects in $\cA'$. Let $(\bC, \bD) : A \to A'$ be a morphism. We
  define a functor $\dQ(\bC, \bD) : \dQ(A) \to \dQ(A')$ as follows. On an
  object $q$ of $\dQ(A)$, define $\dQ(\bC, \bD) (q) = \bD q$. For a morphism
  $\psi : q \to q'$ in $\dQ(A)$, we define $\dQ(\bC, \bD)(\psi) = \bD \psi$.
  \label{dfn:scrQforG}
\end{definition}

Being a full subcategory of a DG category, $\dQ \big(
\objAline{\dC}{\dD}{}{} \big)$ is a DG category.

\begin{lemma}
  $\dQ$ defined above is a (1-)functor from $\cA'$ to $\dgcat$.
  \label{lem:dqisfunc}
\end{lemma}

\begin{proof}
  For reference we have the following diagram
  \begin{equation*}
    \xymatrix{
      \dC \ar[r]^{\iota} \ar[d]^{\bC} & \perf{\dC} \ar[d]^{\perf{\bC}}
      \ar@^{->}@<1pt>[r]^{\bF} & \dD \ar[d]^{\bD} \ar@^{->}@<1pt>[l]^{\bG} &
      \dQ(A) \ar@{_{(}->}[l] \ar[d]^{\dQ(\bC, \bD) = \bD|_{\dQ(A)}} \\
      \dC' \ar[r]^{\iota'} & \perf{(\dC')} \ar@^{->}@<1pt>[r]^{\bF'} & \dD'
      \ar@^{->}@<1pt>[l]^{\bG'} & \dQ(A'). \ar@{_{(}->}[l] \\
    }
  \end{equation*}
  $\bD q$ belongs to $\dQ(A')$ since
  \begin{align*}
    \bG' \bD q &= \perf{\bC} \bG q = \perf{\bC} \iota (c) \text{ for some
    object } c, \text{ by definition of } \dQ(A); \\
    &= \iota' \bC c.
  \end{align*}
  Since $\dQ(A)$ and $\dQ(A')$ are full subcategories, $\dQ(\bC, \bD)$ is a
  functor.
\end{proof}

\begin{definition}
  Let $\cA$ be the full subcategory of $\cA'$ consisting of all objects $A =
  \objAline{\dC}{\dD}{\bF}{\bG}$ such that $\bF|_{\iota\dC}$ takes every
  object of $\iota\dC$ to an object in $\dQ(\cA)$.
\end{definition}

Thus, for $A = \objAline{\dC}{\dD}{\bF}{\bG}$, the adjunction
$\adjunct{\bF}{\perf{\dC}}{\dD}{\bG}$ restricts to an adjunction
$\adjunct{\bF}{\dC}{\dQ(A)}{\bG}$.

To make the exposition clearer, we introduce another category, which we call
$\cB$ as follows.

\begin{definition}
  The objects of $\cB'$ are 4-tuples $(\dC, \dD, \bF, \bG)$ where $\dC$ and
  $\dD$ are DG categories and $\!\adjunct{\bF}{\perf{\dC}}{\dD}{\bG}$ is an
  adjunction. A morphism $\Phi : (\dC, \dD, \bF, \bG) \to (\dC', \dD', \bF',
  \bG')$ in $\cB'$ is a pair of functors $\bC : \dC \to \dC'$ and $\bD : \dD
  \to \dD'$, such that the following diagrams commute
  \begin{equation*}
    \xymatrix{
      \perf{\dC} \ar[r]^{\perf{\bC}} \ar[d]_{\bF} &
      \perf{(\dC')} \ar[d]_{\bF'} &
      &
      &
      \perf{\dC} \ar[r]^{\perf{\bC}} &
      \perf{(\dC')} \\
      \dD \ar[r]^{\bD} &
      \dD' &
      &
      &
      \dD \ar[r]^{\bD} \ar[u]_{\bG} & \dD'. \ar[u]_{\bG'}
    }
  \end{equation*}
  We have an obvious functor $\Xi : \cB' \to \cA'$ such that
  \begin{equation*}
  \Xi((\dC, \dD, \bF, \bG)) = \objA{\dC}{\dD}{\bF}{\bG} \text{ and }
  \Xi((\bC, \bD)) = (\bC, \bD).
  \end{equation*}
  Let $\cB$ be the full subcategory of $\cB'$ consisting of tuples
  $B = (\dC, \dD, \bF, \bG)$ such that $\Xi(B)$ is an object of $\cA$.
  \label{dfn:elagincB}
\end{definition}

We shall construct a natural transformation $\Gamma$ between two functors.
The first functor $\corrmod : \cB \to \Cat$ is constructed as follows.

\tododone{
  Need to modify the construction slightly so that we get the module
  categories in $H^0(\dC)$ instead of $H^0(\perf{\dC})$. Put conditions on
  $\bF$ so that $\bM$ becomes a monad on $\dC$ which factors via $\dQ(\bC,
  \bD)$. Then rewrite all definitions.
}

\begin{definition}
  Given an object $T = (\dC, \dD, \bF, \bG)$ in $\cB$, let $\bM = \bG \circ
  \bF : \dC \to \dC$, where $\bF$ and $\bG$ fit into the restricted
  adjuntion $\adjunct{\bF}{\dC}{\dQ(\Xi(T))}{\bG}$. The Eilenberg-Moore
  construction determines a monad $(\bM, \dmu, \deta)$ on $\dC$.  Applying
  the $2$-functor $H^0$, we get a monad $(M, \mu, \eta)$ on $\cC =
  H^0(\dC)$. We define $\corrmod(T)$ to be the category $\module{M}$ of $M$
  modules in $\cC$.
  
  Now consider a morphism $\Phi : T \to T'$ in $\cB$. Let $T = (\dC, \dD,
  \bF, \bG)$ and $T' = (\dC', \dD', \bF', \bG')$. By definition of $\Phi$,
  we have a map $\bC : \dC \to \dC'$. We define $\corrmod(\Phi)$ to be
  \begin{equation*}
    \corrmod(\Phi)(x, \lambda) = (Cx, C \lambda)
  \end{equation*}
  where $C = H^0(\bC)$ is the induced functor from $\cC = H^0(\dC)$
  to $\cC' = H^0(\dC')$. It is easy to check that $\corrmod(\Phi)$
  is a functor from $\module{M\!}$ to $\module{M'\!}$, where $M = H^0(\bG
  \circ \bF)$ and $M' = H^0(\bG' \circ \bF')$.
\end{definition}

The following is evident.

\begin{lemma}
  $\corrmod : \cB \to \Cat$ is a functor.
\end{lemma}

We define another functor $\tilde{Q}$.
\begin{definition}
  Let $\tilde{Q} = \dQ \circ \Xi : \cB \to \Cat$ be the functor, where $\dQ
  : \cA \to \dgcat$ was defined in \ref{dfn:scrQforG} and $\Xi : \cB \to
  \cA$ was defined in \ref{dfn:elagincB}.
\end{definition}

We define an association $\Gamma$ as follows

\begin{definition}
  For an object $T = (\dC, \dD, \bF, \bG)$ of $\cB$, define $\Gamma T :
  H^0(\tilde{Q}(T)) \to \corrmod(T)$ to be the functor defined as follows.
  Note that, since $T$ is an object of $\cB$, we have an adjunction
  \begin{equation*}
    \adjunct{\bF}{\dC}{\tilde{Q}(T)}{\bG}.
  \end{equation*}
  Taking $H^0$, and setting $M = H^0(\bG \circ \bF)$, we get a diagram
  \begin{equation*}
    \xymatrix{
      \cC = H^0(\dC) \ar@<1pt>@^{->}[d] \ar@<1pt>@^{->}[rrd] &
      &
      \\
      H^0(\tilde{Q}(T)) \ar@<1pt>@^{->}[u] \ar@{-->}[rr]^{\Gamma} &
      &
      \module{M} \ar@<1pt>@^{->}[llu]
    }
  \end{equation*}
  where $\Gamma$ is the universal functor coming from lemma
  \ref{lem:emenrmnd}.

\end{definition}

\begin{proposition}
  $\Gamma : H^0(\tilde{Q}) \to \corrmod$ is a natural transformation.
\end{proposition}
\begin{proof}
  Given the definitions, this is a routine computation.
\end{proof}

\subsection{Elagin's construction} \label{ssc:elagincst}

\begin{theorem}
  Suppose $\cC$ is a triangulated category and $(M, \mu, \eta)$ is an exact
  monad on $\cC$. Let $\idemcomp{\cC}$ be the idempotent completion of $\cC$
  and $(\idemcomp{M}, \idemcomp{\mu}, \idemcomp{\eta})$ be the monad induced
  on $\idemcomp{\cC}$ (see definition \ref{dfn:idmcmplt}). Assume the
  following.
  \begin{enumerate}
  \item There exists an enhancement $(\dC, \varepsilon)$ of $\cC$ and a
    lift $(\bM, \dmu, \deta)$ of $M$ on $\dC$. Let $\perf{\dC}$ be the
    category of perfect DG $\dC$-modules (see \ref{dfn:perfdgct}) and let
    $(\perf{\bM}, \perf{\dmu}, \perf{\deta})$ be the induced DG monad on
    $\perf{\dC}$.
    \item $\module{\idemcomp{M}}$ has an enhancement $(\dD, \delta)$ which
      admits a pair of adjoint DG functors
      $\adjunct{\bF}{\perf{\dC}}{\dD}{\bG}$ satisfying $\perf{\bM} = \bG
      \circ \bF$ and $H^0(\bG) = G_{\idemcomp{M}}$ (see \ref{dfn:emcons})
      such that in the diagram
      \begin{equation*}
        \xymatrix{
          H^0(\dD) \ar[rr]^{\delta} \ar@<1pt>@^{->}[dr]^{H^0(\bG)} &
          &
          \module{\idemcomp{M}} \ar@<1pt>@^{->}[dl]^{G_{\idemcomp{M}}} \\
          &
          \idemcomp{\cC} \ar@<1pt>@^{->}[ul]^{H^0(\bF)}
          \ar@<1pt>@^{->}[ur]^{F_{\idemcomp{M}}} &
        }
      \end{equation*}
      $F_{\idemcomp{M}} = \delta \circ H^0(\bF)$ and $H^0(\bG) =
      G_{\idemcomp{M}} \circ \delta$.
  \item We also assume that $(\dC, \dD, \bF, \bG)$ is an object of
    $\cB$ (see definition \ref{dfn:elagincB}).
  \end{enumerate}
  Under the above hypothesis, if $B = (\dC, \dD, \bF, \bG)$, the following
  holds
  \begin{enumerate}
    \item $\Gamma(B) : H^0(\tilde{Q}(B)) \to \corrmod(B) \cong
      \module{M}$ is an equivalence of categories.
    \item $\module{M}$ admits a triangulated structure with respect to which
    $F_M$ and $G_M$ (see \ref{dfn:emcons}) becomes exact.
    \item Furthermore, $\module{M}$ admits an enhancement.
  \end{enumerate}
  \label{thm:Mmodenhd}
\end{theorem}

\begin{proof}
  Let $F = H^0(\bF)$ and $G = H^0(\bG)$. If $A =
  \delta(H^0(\tilde{Q}(B))$, $A$ would fit into the following two diagrams
  \begin{equation*}
    \xymatrix@=.3cm{
      &
      \idemcomp{\cC} \ar[dl]_{F} \ar[dr]^{F_{\idemcomp{M}}} &
      &
      &
      &
      \idemcomp{\cC} &
      \\
      H^0(\dD) \ar[rr]^(.4){\delta} &
      &
      \module{\idemcomp{M}} &
      &
      H^0(\dD) \ar[rr]^(.4){\delta} \ar[ru]^{G} &
      &
      \module{\idemcomp{M}} \ar[lu]_{G_{\idemcomp{M}}} \\
      &
      &
      &
      &
      &
      &
      \\
      &
      \cC \ar[dl]_{F} \ar[dr]^{F_{\idemcomp{M}}} 
      \ar@{>->}[uuu]|(.64){\hole} &
      &
      &
      &
      \cC \ar@{>->}[uuu]|(.64){\hole} &
      \\
      H^0(\tilde{Q}(B)) \ar[rr]^{\delta} \ar@{>->}[uuu] &
      &
      A \ar@{>->}[uuu] &
      &
      H^0(\tilde{Q}(B)) \ar[rr]^{\delta} \ar[ru]^{G} \ar@{>->}[uuu] &
      &
      A \ar[lu]_{G_{\idemcomp{M}}} \ar@{>->}[uuu] 
    }
  \end{equation*}
  Since, $G_{\idemcomp{M}}$ maps $A$ in to $\cC$, it is easy to see that $A
  \cong \module{M}$; and since $\delta$ is an equivalence, universal
  property of $\Gamma(B)$ implies that $\Gamma(B) = \delta :
  H^0(\tilde{Q}(B)) \to \module{M} = \corrmod(B)$ is an equivalence; thus
  proving the first part.

  We know $\module{M}$ is equivalent to $H^0(\tilde{Q}(B))$. In order to
  show that $\module{M}$ is triangulated, we show that $\tilde{Q}(B)$ is
  pretriangulated. Recall that both $\dC$ and $\dD$ are pre-triangulated.
  Following remark \ref{rmk:altdscpt}, we prove that shifts of objects in
  $\tilde{Q}(B)$ and cones of degree zero closed morphisms in $\tilde{Q}(B)$
  are homotopically equivalent to some object in $\tilde{Q}(B)$.

  The dg functor $\bG : \tilde{Q}(B) \to \dC$ induces a functor $\pretr{\bG}
  : \pretr{(\tilde{Q}(B))} \to \pretr{\dC}$. Let $d$ be an object of
  $\tilde{Q}(B)$ and $n \in \Z$. Now $\pretr{\bG}(d[n]) =
  (\pretr{\bG}(d))[n]$ is homotopy equivalent to some $c_n$ in $\dC$, since
  $\dC$ is pretriangulated. Thus $d[n]$ belongs to $\tilde{Q}(B)$, by
  definition.

  By the construction of cone (see \ref{dfn:twistdcp}), it is
  clear that $\pretr{\bG}(\Cone(f)) = \Cone (\pretr{\bG}(f))$. The degree
  zero closed map $\pretr{\bG}(f) = \bG f$ is a morphism in $\dC$, which is
  pretriangulated. Thus $\Cone (\pretr{\bG}(f))$ is homotopy equivalent to
  some object $c$ in $\dC$. In other words, $\pretr{\bG}(\Cone(f))$ is
  homotopy equivalent to some object in $\dC$. Therefore, $\Cone(f)$
  belongs to $\tilde{Q}(B)$. Thus $\tilde{Q}(B)$ is pretriangulated.

  Thus $\module{M} = \corrmod(B)$, which is equivalent to
  $H^0(\tilde{Q}(B))$ is triangulated.  Furthermore, the enhancement is
  given by $(\tilde{Q}(B), \Gamma(B))$. This proves 2 and 3.
\end{proof}

As an application we show how the above theorem easily gives us corollary
6.10, Theorem 8.9 and Corollary 8.10 in Elagin's paper
\cite{elagin:onequivtrcat}. Now we recall the definition of action of a
group on a category.

\begin{definition}
  Suppose $G$ is a finite group of order $\abs{G}$. Let $\cV$ be a
  $\Z[1/\abs{G}]$-linear symmetric, monoidal category. For a $\cV$-enriched
  category $\cC$, an action of $G$ on $\cC$ is defined to be a family of
  $\cV$-autoequivalences $\set{\varphi_g}{\varphi_g : \cC \to \cC, g \in
  G}$, such that we have a compatible family of $\cV$-natural isomorphisms
  $\eta_{g, h} : \varphi_g \circ \varphi_h \xrightarrow{\cong}
  \varphi_{hg}$, in the sense for all $f, g, h \in G$, the diagram
  \begin{equation*}
    \xymatrix{
      \varphi_f \varphi_g \varphi_h \ar[r]^{\varphi_f
        \eta_{g, h}} \ar[d]_{\eta_{f, g} \varphi_h} & \varphi_f
        \varphi_{hg} \ar[d]^{\eta_{f, hg}} \\
        \varphi_{gf} \varphi_h \ar[r]^{\eta_{gf, h}} & \varphi_{hgf}
    }
  \end{equation*}
  commutes.
  \label{dfn:gactonct}
\end{definition}

\begin{definition}
  For $G$ a finite group, $\cV$ a $\Z[1/\abs{G}]$-linear category and for
  $\cC$ a $\cV$-enriched category with an action of $G$ as in definition
  \ref{dfn:gactonct}, one defines the $\cV$-enriched category $\cC^G$ of
  $G$-equivariant objects as follows.  An object of $\cC^G$ is a tuple
  consisting of
  \begin{enumerate}
    \item an object $A$ of the category $\cC$,
    \item a collection of isomorphisms in $\cC$, $\set{\lambda_g : A \to
      \varphi_g A}{g \in G}$
  \end{enumerate}
  such that for all $g$ and $h$ in $G$, the diagram
  \begin{equation*}
    \xymatrix{
      A \ar[r]^{\lambda_g} \ar[d]_{\lambda_{hg}} & \varphi_g A
      \ar[d]^{\varphi_g(\lambda_h)} \\
      \varphi_{hg} A & \varphi_g \varphi_h A \ar[l]^{\eta_{g, h}}_{\cong}
    }
  \end{equation*}
  commutes.

  A morphism $f : (A, \setlist{\lambda_g : A \to \varphi_g A}) \to (B,
  \setlist{\tau_g : B \to \varphi_g B})$, between two objects of $\cC^G$, is
  a morphism $f : A \to B$ in $\cC$ which is compatible with the action of
  $G$. In other words:
  \begin{equation*}
    \xymatrix{
      A \ar[r]^{\lambda_g} \ar[d]_{f} & \varphi_g A \ar[d]^{\varphi_g f} \\
      B \ar[r]^{\tau_g} & \varphi_g B
    }
  \end{equation*}
  commutes for all $g \in G$.
  \label{dfn:gequiobj}
\end{definition}

Elagin \cite[example 3.7]{elagin:onequivtrcat} defines the following two
functors.

\begin{definition}
  $p^* : \cC^G \to \cC$ is the forgetful functor $p^*(A,
  \setlist{\lambda_g}) = A$. If $f : (A, \setlist{\lambda_g}) \to (B,
  \setlist{\tau_g})$ is a morphism in $\cC^G$, then $p^*(f)$ is the
  underlying morphism $f : A \to B$.

  $p_* : \cC \to \cC^G$ is defined by
  \begin{align*}
    p_*(A) &= \left( \bigoplus_{g \in G} \varphi_g A, \setlist{\xi_h}
    \right) \\
    &\text{ where } \xi_h = \oplus_{g} \eta_{h, g}^{-1} : \oplus_{g}
    \varphi_{gh} A \rightarrow \oplus_{g} \varphi_h (\varphi_g A) =
    \varphi_h (\oplus_{g} \varphi_g A); \text{ and } \\
    p_*(f : A \to B) &= \bigoplus_{g \in G} \varphi_g f.
  \end{align*}
  \label{dfn:gequivcC}
\end{definition}

\begin{lemma}
  $p^*$ is both right and left adjoint to $p_*$.
  \label{lem:geqvpadj}
\end{lemma}
The above lemma is proved in \cite[lemma 3.8]{elagin:onequivtrcat}.

\begin{definition}
  Given $G$, $\cV$, $\cC$ and $\cC^G$ as in definitions \ref{dfn:gactonct}
  and \ref{dfn:gequivcC}, let $M_G : \cC \to \cC$ be the functor $M_G = p^*
  \circ p_*$.
  \label{dfn:gequivMG}
\end{definition}

Now let us restrict to the case when $\cV$ is the category of dg-modules
over some field $k$. We get the following result easily.

\begin{definition}
  Let $G$ be a finite group of order $\abs{G}$ and suppose $\cC$ is a
  triangulated category over an $\Z[1/\abs{G}]$-algebra $k$. Assume $G$ acts
  on $\cC$. A $G$-equivariant enhancement $(\dC, \varepsilon)$ is a DG
  category $\dC$ which admits a $G$ action via DG autoequivalences, along
  with a $G$-equivariant equivalence $\varepsilon : H^0(\dC) \cong \cC$. As
  in \cite[section 8]{elagin:onequivtrcat}, we assume that the natural
  isomorphisms $\eta_{g,h}$ for $\dC$ (cf.~\ref{dfn:gactonct}) are closed of
  degree $0$. Also, the $G$-equivariant objects in $\dC^G$ are defined as
  pairs $(A, \set{\dlambda_g}{g \in G})$ as in \ref{dfn:gequiobj}, with the
  additional condition that $\dlambda_g : A \to \varphi_g A$ are all closed
  morphisms of degree $0$.
\end{definition}

\begin{theorem}
  Suppose $G$ is a finite group of order $\abs{G}$ and $\cC$ is a
  triangulated category over an $\Z[1/\abs{G}]$-algebra $k$. We assume that
  $G$ acts on $\cC$. If $\cC$ admits an $G$-equivariant enhancement $(\dC,
  \varepsilon)$, then $\cC^G$ is a triangulated category, and admits an
  enhancement.
  \label{thm:geqvtcat}
\end{theorem}

\begin{proof}
  \tododone{
  This follows from the fact that $M_G$ on $\cC$ admits a lift $\bM_G$
  as a monad and hence theorem \ref{thm:Mmodenhd}(3) is applicable.
  }
  The construction of $p_*$ and $p^*$ in \ref{dfn:gequivcC} applied for
  the case where $\cV$ is the category of DG $k$-modules gives us two
  DG functors $\boldp_* : \dC \to \dC^G$ and $\boldp^* : \dC^G \to \dC$. Let
  $\bM = \boldp^* \circ \boldp_*$. It is clear that $H^0(\bM) =
  H^0(\boldp^*) \circ H^0(\boldp_*) = p^* p_* = M$. Thus $M$ admits a lift
  $\bM$; and hence condition 1 of theorem \ref{thm:Mmodenhd} is satisfied.

  We refer to \cite[lemma 8.6 and theorem 8.7]{elagin:onequivtrcat} for the
  proof that $\module{\idemcomp{M}} \cong (\idemcomp{\cC})^G \cong
  H^0(\moduledg{\perf{\bM}}) \cong H^0((\perf{\dC})^G) \cong
  H^0(\perf{(\dC^G)})$. This verifies condition 2 of theorem
  \ref{thm:Mmodenhd}.

  Note that in our case, $\bF = \boldp_*$ and $\bG = \boldp^*$. Let
  $B = (\dC, \perf{(\dC^G)}, \boldp_*, \boldp^*)$. For an object $c$ in
  $\dC$, $\bG \bF c = \boldp^* \boldp_* c = \oplus \varphi_g c$ is actually
  an object in $\dC$. Thus, condition 3 of theorem \ref{thm:Mmodenhd} is
  also satisfied.

  Thus, by the theorem $\module{M}$ is triangulated and admits an enhancement
  $(\tilde{Q}(B), \Gamma(B))$.
\end{proof}

\begin{remark}
  Elagin \cite[theorem 6.9 and corollary 6.10]{elagin:onequivtrcat} proves
  that $\module{M} = \cC^G$ is triangulated even when $\cC$ has an
  enhancement, which may not be $G$-equivariant. He uses separability of
  $M$ and the fact that a triangulated category having a DG enhancement is
  triangulated of order $N$ for any $N > 0$.
\end{remark}

In particular, we have the following geometric results.

\begin{corollary}
  Suppose $X$ is quasi-projective smooth variety over a field $k$ and $G$ is
  a finite group of order $\abs{G}$ such that $(\charac k, \abs{G}) = 1$.
  \begin{enumerate}
    \item Assume that $G$ acts on $X$.  Let $\cD^b_G(X)$ denote the bounded
      derived category of $G$-equivariant coherent sheaves. Then with the
      induced action of $G$ on $\cD^b(X)$, we have
      \begin{equation*}
        \cD_G^b(X) \cong \cD^b(X)^G.
      \end{equation*}
    \item Suppose $G$ is a finite subgroup of $\Pic(X)$. Then $G$ acts on
      $\coh{X}$ by tensoring into line bundles in $G$. If $Y$ is the
      relative spectrum
      \begin{equation*}
        Y = \Spec_X \left(\bigoplus_{L \in G} L^{-1}\right),
      \end{equation*}
      then $\coh{X}^G \cong \coh{Y}$ and hence $\cD^b_G(X) \cong
      \cD^b(Y)$.
  \end{enumerate}
\end{corollary}

\begin{proof}
  Since $X$ has enough locally free sheaves,
  \begin{equation*}
  \cD^b(X) \cong \Dperf(X) = \cK^b(\lffr(X)).
  \end{equation*}
  where $\lffr(X)$ denotes the additive category of locally free, finite
  rank sheaves over $X$ and $\cK^b(\underline{\ \ })$ denotes the homotopy
  category of chain complexes. If $\lffr_G(X)$ denotes the category of
  locally free, finite rank $G$-equivariant sheaves on $X$,
  \begin{equation*}
  \cD^b_G(X) \cong \cK^b(\lffr_G(X))
  \end{equation*}
  Thus, $\cD^b(X)$ is a triangulated category admitting a $G$-equivariant
  enhancement $\dC^b(\lffr(X))$, where $\dC^b(\lffr(X))$ is the DG category
  of bounded chain complexes of objects from $\lffr(X)$. Theorem
  \ref{thm:geqvtcat} implies that $\cD^b(X)^G$ is triangulated and admits an
  enhancement.
  
  The construction of $p_*$ and $p^*$ in definition \ref{dfn:gequivcC} gives
  us a monad $M = p^* \circ p_* : \lffr(X) \to \lffr(X)$. Now we are in the
  situation of example \ref{exl:monadabc}. Let $\bM_{dg}$ be the monad
  induced on $\dC^b(\lffr(X))$ by $M$. By lemma \ref{lem:egcomplx},
  $H^0(\moduledg{\bM_{dg}}) \cong \cK^b(\module{M})$. Now $\module{M} =
  \lffr_G(X)$ and $\moduledg{\bM_{dg}} \cong \left( \dC^b(\lffr(X))
  \right)^G$. Thus,
  \begin{equation*}
  \cD^b(X)^G = H^0(\dC^b(\lffr(X))^G) \cong \cK^b(\lffr_G(X)) \cong
  \cD^b_G(X).
  \end{equation*}

  The second part follows similarly. \tododone{Check we don't need extra
  conditions on the action of $G$ on $X$.}
\end{proof}

\section{Twisted derived categories} \label{sec:twistedc}

\tododone{Define twisted derived category}
We recall the definition of twisted sheaves and twisted derived categories.
References are \cite[section 4]{caldararu:twistdercatell3fld} and his
thesis \cite{caldararu:thesis}.
\begin{definition}
  Let $X$ be a scheme. Suppose $\alpha \in \cech^2(X, \sheafO_X^*)$ is given
  along an open cover $\setlist{U_i}_{i \in I}$ by sections $\alpha_{i,j,k}
  \in \Gamma(U_i \cap U_j \cap U_k, \sheafO_X^*)$. An $\alpha$-twisted sheaf
  $\sheafF$ consists of a pair $(\setlist{\sheafF_i}_{i \in I},
  \setlist{\varphi_{i,j}}_{i, j \in I})$, where $\sheafF_i$ are sheaves on
  $U_i$ and $\varphi_{i, j} : \sheafF_j|_{U_i \cap U_j} \to
  \sheafF_i|_{U_i \cap U_j}$ are isomorphisms such that
  \begin{enumerate}
  \item $\varphi_{i, i} = \id_{\sheafF_i}$,
  \item $\varphi_{i, j} = \varphi_{j, i}^{-1}$, and
  \item $\varphi_{i, j} \circ \varphi_{j, k} \circ \varphi_{k, i} =
    \alpha_{i, j, k} \id_{\sheafF_i|_{U_i \cap U_j \cap U_k}}$.
  \end{enumerate}
  Morphism of twisted sheaves is defined in an obvious manner: they are a
  collection of morphisms of sheaves on $U_i$ which are compatible with the
  twisting.  Given $X$ and $\alpha$, the twisted sheaves form an abelian
  category, which we denote by $\Mod(X, \alpha)$. The subcategory of
  $\Mod(X, \alpha)$ consisting of twisted sheaves with all the $\sheafF_i$'s
  coherent is the category of coherent twisted sheaves, which is denoted by
  $\coh{X, \alpha}$. $\coh{X, \alpha}$ is an abelian category. The category of 
  locally free, finite rank twisted sheaves will be denoted by
  $\lffr(X, \alpha)$.

  Let $\cD^b(X, \alpha)$ denote the bounded derived category of complexes of
  $\alpha$-twisted sheaves with coherent cohomology.
\end{definition}

\begin{remark}
  The categories $\Mod(X, \alpha)$ and $\coh{X, \alpha}$ do not depend on
  the choice of open cover or on the choice of the cocycle
  $\setlist{\alpha_{i,j,k}}$. Remark
  \cite[4.5]{caldararu:twistdercatell3fld} and \cite[lemma
  2.1.4]{caldararu:thesis} implies that if $\alpha \in \Br(X)$ and $\coh{X}$
  has enough locally free sheaves of finite rank, then  $\coh{X, \alpha}$
  has enough locally free $\alpha$-twisted sheaves of finite rank. Thus, in
  this case,
  \begin{equation*}
  \cD^b(X, \alpha) \cong \cK^{-, b}(\lffr(X, \alpha))
  \end{equation*}
  where $\cK^{-, b}$ is the homotopy category of bounded above complexes
  which are exact at all but finitely many degrees.
  \label{rmk:tderKtwi}
\end{remark}

\tododone{Refer to the result $\cD_b(X, \alpha) \cong \cK^-(\lffr_{\alpha}(X))$}

We need a few results before we prove the main theorem of this section.

\begin{proposition}
  A sheaf of Azumaya algebras $\cA$ over $\sheafO_X$ on a scheme $X$ is
  separable over $\sheafO_X$.
  \label{prp:azumasep}
\end{proposition}

\begin{proof}
  This follows from section 3 of chapter 2 in \cite{di:sepalgcommrings}.
\end{proof}

\begin{proposition}
  If $S$ is a separable $R$ algebra, the forgetful functor $\text{ff} :
  \module{S} \to \module{R}$ is separable.
  \label{prp:sepalgff}
\end{proposition}
\begin{proof}
  This follows from proposition 1.3 in \cite{nvo:sepfuncgrring}.
\end{proof}

\begin{proposition}
  If $\cA$ is an Azumaya algebra over $X$ corresponding to $\alpha \in
  \cech^2(X, \sheafO_X^*)$, then the categories $\Mod(X, \alpha)$ and
  $\Mod(\cA)$ are equivalent.
  \label{prp:ashftshf}
\end{proposition}
\begin{proof}
  This is Theorem 1.3.7 in \cite{caldararu:thesis}.
\end{proof}

\begin{lemma}
  Let $\alpha$ be an element of the Brauer group $\Br(X)$ corresponding to
  the Azumaya algebra $\cA_{\alpha}$. Then 
  \begin{equation*}
    \cD^b(X, \alpha) \cong \module{M_{\alpha}},
  \end{equation*}
  where $M_{\alpha}$ is the monad $M_{\alpha} (\cF) = \cF \otimes
  \cA_{\alpha}$ on $\cD^b(X)$.
  \label{lem:twisderc}
\end{lemma}
\begin{proof}

  This can be directly proved using theorem \ref{thm:Mmodntrg}(4) as
  follows. Consider the monad $M_{\alpha}$ on $\cD^b(X)$ defined by
  \begin{equation*}
  M_{\alpha}(\cF) = \cF \otimes_{\sheafO_X} \cA_{\alpha}
  \end{equation*}
  where $\otimes$ denotes the left derived tensor. We have an adjunction
  \begin{equation*}
  \adjunct{F}{\cD^b(X)}{\cD^b(\module{\cA_{\alpha}})}{G}
  \end{equation*}
  where
  \begin{equation*}
  F(\cF) = \cF \otimes_{\sheafO_X} \cA_{\alpha} \qquad G(\cG) =
  \iota^* \cG = \Hom_{\cA_{\alpha}}(\cA_{\alpha}, \cG),
  \end{equation*}
  where $\iota : \sheafO_X \to \cA_{\alpha}$ is the algebra map and
  $\iota^*$ is the restriction of scalars functor. Also $M_{\alpha} = G
  \circ F$.

  Since Azumaya algebras are separable by \ref{prp:azumasep}, $G$ is
  separable by \ref{prp:sepalgff}. Thus by Balmer's theorem
  \ref{thm:Mmodntrg}(4) we get that $\module{M_{\alpha}} \cong
  \cD^b(\module{\cA_{\alpha}}) \cong \cD^b(X, \alpha)$. The last equivalence
  follows from \ref{prp:ashftshf}.
\end{proof}

\begin{theorem}
  Suppose $\cA$ is an Azumaya algebra over $X$ associated to $\alpha \in
  \cech^2(X, \sheafO_X^*)$. Then $\cD^b(X, \alpha)$ admits an enhancement.
  \label{thm:twistenh}
\end{theorem}

\begin{proof}
  Let us denote the functor $\cF \mapsto \cF \otimes \cA_{\alpha}$ from
  $\lffr(X)$ to itself by $M_{\alpha}$. It is clear that since
  $\cA_{\alpha}$ is an algebra, $M_{\alpha}$ is a monad. Consider
  the induced functor
  \begin{equation*}
  \bM_{\alpha} : \cC_{dg}^{-, b} (\lffr(X)) \to \cC_{dg}^{-, b}
  (\lffr(X)).
  \end{equation*}
  Note $\module{M_{\alpha}} \cong \lffr(X, \alpha)$. Thus
  $H^0(\moduledg{\bM_{\alpha}}) \cong \cK^{-, b}(\lffr(X, \alpha)) \cong
  \cD^b(X, \alpha)$, where the first equivalence comes from lemma
  \ref{lem:egcomplx} and the second one comes from remark
  \ref{rmk:tderKtwi}. This demonstrates an enhancement for $\cD^b(X,
  \alpha)$.
\end{proof}

\section{Interaction of two monad actions} \label{sec:compatmnd}

In this section, we explore conditions on two monads which ensure that their
composition is also a monad.

\begin{definition}
  Let $\cC$ be a category. Two monads $(M_1, \mu_1, \eta_1)$ and $(M_2,
  \mu_2, \eta_2)$ are said to be compatible if the following holds.
  \begin{enumerate}
  \item There exists a natural isomorphism $\sigma : M_1 M_2 \to M_2 M_1$.
  \item It is compatible with $\mu_1$ and $\mu_2$ as is shown by the
    commutativity of 
    \begin{equation*}
    \xymatrix{
      M_1^3 M_2^3 \ar[rr]^{(\mu_1 M_1) \ast (\mu_2 M_2)}
      \ar[d]_{M_1^2 \sigma M_2^2} &
      &
      M_1^2 M_2^2 \ar[dd]^{M_1 \sigma M_2} \\
      M_1^2 M_2 M_1 M_2^2 \ar[d]_{M_1^2 M_2 \sigma M_2} &
      &
      \\
      M_1^2 M_2^2 M_1 M_2 \ar[rr]^{(\mu_1 \ast \mu_2) M_1 M_2} &
      &
      M_1 M_2 M_1 M_2
    }
    \end{equation*}
    along with a similar diagram with roles of $M_1$ and $M_2$ reversed.
  \item Also the following are commutative
    \begin{equation*}
    \xymatrix{
      M_1^2 M_2 \ar[r]^{\mu_1 M_2} \ar[d]_{M_1 \sigma} &
      M_1 M_2 \ar[dd]^{\sigma} &
      M_2^2 M_1 \ar[r]^{\mu_2 M_1} &
      M_2 M_1 &
      \\
      M_1 M_2 M_1 \ar[d]_{\sigma M_1} &
      &
      M_2 M_1 M_2 \ar[u]^{M_2 \sigma} &
      &
      \\
      M_2 M_1^2 \ar[r]^{M_2 \mu_1} &
      M_2 M_1 &
      M_1 M_2^2 \ar[r]^{M_1 \mu_2} \ar[u]^{\sigma M_2} &
      M_1 M_2 \ar[uu]_{\sigma} &
    }
    \end{equation*}
    along with two similar diagrams with $M_1 M_2^2$ and $M_2 M_1^2$ on the
    top left corner.
  \item We have two similar conditions for $\eta$. The first one being
    compatibility with horizontal composition.
    \begin{equation*}
      \xymatrix{
        M_1 M_2 \ar[rr]^{M (\eta_1 \ast \eta_2)} \ar[drr]^{(M_1 \eta_1) \ast
        (M_2 \eta_2)} \ar[d]_{(\eta_1 \ast \eta_2) M} &
        &
        M_1 M_2 M_1 M_2 \\
        M_1 M_2 M_1 M_2 &
        & 
        M_1^2 M_2^2 \ar[ll]^{M_1 \sigma M_2} \ar[u]_{M_1 \sigma M_2}.
      }
    \end{equation*}
  \item We also need the following compatibility.
    \begin{equation*}
      \xymatrix{
        M_2 \ar[r]^-{\eta_1 M_2} \ar[d]_{M_2 \eta_1} &
        M_1 M_2 \ar[dl]_{\sigma} \\
        M_2 M_1 &
        M_1 \ar[u]_{M_1 \eta_2} \ar[l]_-{\eta_2 M_1}
    }
    \end{equation*}
  \end{enumerate}
\end{definition}

\begin{remark}
  If we consider the map $\sigma_z : M_1 M_2 z \to M_2 M_1 z$, and use
  the definition of natural transformation on $\sigma$, we get the following
  diagram (where we omit the z).
  \begin{equation*}
    \xymatrix{
      M_1 M_2 M_1 M_2 \ar[rr]^{\sigma M_1 M_2} \ar[d]_{M_1 M_2 \sigma} &
      &
      M_2 M_1^2 M_2 \ar[d]^{M_2 M_1 \sigma}\\
      M_1 M_2^2 M_1 \ar[rr]^{\sigma M_2 M_1} &
      &
      M_2 M_1 M_2 M_1
    }
  \end{equation*}
\end{remark}

\tododone{The above conditions are enough to conclude the first two parts of
the following proposition. Do they imply the third?}

We need the following proposition.

\begin{proposition}
  Let $\dC$ be a DG category with two compatible monads $(\bM_1, \dmu_1,
  \deta_1)$ and $(\bM_2, \dmu_2, \deta_2)$.
  \begin{enumerate}
    \item $(\bM_1 \bM_2, (\dmu_1 \ast \dmu_2) \circ (\bM_1 \sigma^{-1}
      \bM_2), \deta_1 \ast \deta_2)$ is a monad on $\cC$. By symmetry, so is
      $(\bM_2 \bM_1, (\dmu_2 \ast \dmu_1) \circ (\bM_2 \sigma^{-1} \bM_1),
      \deta_1 \ast \deta_2)$.
    \item $\bM_2$ (respectively $\bM_1$) induces a monad say
      $\overline{\bM}_2$ (respectively $\overline{\bM}_1$) on
      $\module{\bM}_1$ (respectively $\module{\bM}_2$).
    \item $\module{\overline{\bM}_2} \cong \module{(\bM_1 \bM_2)}$ and
      $\module{\overline{\bM}_1} \cong \module{(\bM_2 \bM_1)}$. Further,
      $\module{(\bM_1 \bM_2)} \cong \module{(\bM_2 \bM_1)}$.
  \end{enumerate}
  \label{prp:compmond}
\end{proposition}

\begin{proof}
  This was done by Beck, for the case of additive categories, in the
  proposition stated in \cite[page 126 -- 127]{beck:distributivelaws}, and
  the work done before that. The same proof can be adapted for any
  $\cV$-category and hence for DG categories. The last part of 9.3(3) is not
  in Beck, but it follows from the fact that $\sigma$ is a natural
  isomorphism.
\end{proof}

As an application, we have the following two results.

\tododone{Taking supports commute with twisting}

\begin{definition}
  Let $\cdgr{\cA}{\cB}$ be the full subcategory of $\cdg{\cA}$ (see example
  \ref{exl:monadabc}) consisting of those complexes, whose homologies lie in
  $\cB$.
  \label{dfn:dgcpxwsp}
\end{definition}

\begin{lemma}
  Suppose $\cA$ is a Grothendieck abelian category and $\cB$ is a Serre
  subcategory. Suppose every object of $\cdg{\cA}$ admits a K-injective
  resolution. Let $\cI$ be the full subcategory of all K-injective objects
  in $\cdgr{\cA}{\cB}$. Then there exists a DG monad $\bM_{d}$ on
  $\cdgr{\cA}{\cB}$ such that $\moduledg{\bM_{d}} \cong \cI$.
  \label{lem:abcatwsp}
\end{lemma}

\begin{proof}
  The inclusion functor $\boldi : \cdgr{\cI}{\cB} \hookrightarrow
  \cdgr{\cA}{\cB}$ has a left adjoint $\bI$ which takes every complex in
  $\cdgr{\cA}{\cB}$ to its K-injective resolution. Let $\bM_d = \boldi \circ
  \bI$. $\dmu$ and $\deta$ are defined as in the Eilenberg-Moore
  construction. This has the required properties, as any object in
  $\moduledg{\bM_d}$, $\lambda : \bI x \to x$ gives a splitting of $x \to
  \bI x$ which proves that $x$ is also K-injective.
\end{proof}

\begin{lemma}
  Let $\cA$ be a Grothendieck, abelian category with tensor and $\cB$ be a
  Serre category. Let $\one$ be the identity for the tensor. Consider an
  object $A^{\bullet}$ in $\cdg{\cA}$, which is K-flat. Suppose
  $A^{\bullet}$ is a ring object in $\cdg{\cA}$, that is, $B^{\bullet}
  \mapsto B^{\bullet} \otimes_{\cdg{\cA}} A^{\bullet}$ is a DG monad on
  $\cdg{\cA}$.
  

  We assume that for any object $B^{\bullet}$ in $\cdgr{\cA}{\cB}$,
  $A^{\bullet} \otimes B^{\bullet}$ is also an object of $\cdgr{\cA}{\cB}$.
  We further assume that $A^{\bullet}$ has a K-flat dual
  $\dual{(A^{\bullet})}$, that is
  \begin{equation*}
	\Hom_{\cdg{\cA}}(M^{\bullet}, A^{\bullet} \otimes N^{\bullet}) \cong
	\Hom_{\cdg{\cA}}(M^{\bullet} \otimes \dual{(A^{\bullet})}, N^{\bullet}).
  \end{equation*}
  Then $\bM_A = A^{\bullet} \otimes \_$ induces a monad structure on
  $\cdgr{\cA}{\cB}$. This monad structure is compatible with $\bM_d$ defined
  above.
  \label{lem:tensring}
\end{lemma}

\begin{proof}
  Every object of $\cdg{\cA}$ admits a K-injective resolution, see for
  example \cite[Theorem 5.4]{tls:loccatres}.

  That $\bM_A$ is a monad, follows from the fact that $A$ is a ring object.
  We only need to show that $\bM_A$ and $\bM_d$ are compatible.

  First we claim that if $I^{\bullet}$ is K-injective, so is $I^\bullet
  \otimes A^{\bullet}$. This is because, for any acyclic $M^{\bullet}$, and
  for any K-injective $I^{\bullet}$,
  \begin{equation*}
	\Hom_{\cdg{\cA}}(M^{\bullet}, I^{\bullet} \otimes A^{\bullet}) =
	\Hom_{\cdg{\cA}}(M^{\bullet} \otimes \dual{(A^{\bullet})}, I^{\bullet})
	= 0
  \end{equation*}
  as by definition of K-flat, $M^{\bullet} \otimes \dual{(A^{\bullet})}$ is
  acyclic.
  
  Suppose for a complex $M^{\bullet}$ in $\cdgr{\cA}{\cB}$,
  $I^{\bullet}(M^{\bullet})$ be its K-injective resolution in $\cdg{\cA}$.
  Being quasi-isomorphic to $M^{\bullet}$, $I^{\bullet}(M^{\bullet})$ also
  belongs to $\cdgr{\cA}{\cB}$. $A^{\bullet}$ being flat,
  $I^{\bullet}(M^{\bullet}) \otimes A^{\bullet}$ is a K-injective resolution
  of $M^{\bullet} \otimes A^{\bullet}$. This gives us a homotopy equivalence
  \begin{equation*}
	I^{\bullet}(M^{\bullet} \otimes A^{\bullet}) \sim
	I^{\bullet}(M^{\bullet}) \otimes A^{\bullet}.
  \end{equation*}
  In other words,
  \begin{equation*}
	\bM_d (\bM_A (M^{\bullet})) \sim \bM_A (\bM_d (M^{\bullet})),
  \end{equation*}
  which is easily checked to be natural, proving the required compatibility.
\end{proof}

\begin{proposition}
  Let $X$ be a quasi-projective variety and let $Z \subset X$ be a closed
  subvariety. Suppose $\cA_{\alpha}$ is an Azumaya algebra over $X$
  associated to $\alpha \in \cech^2(X, \sheafO_X^*)$. Then $\cD_Z(X,
  \alpha)$ admits an enhancement.
\end{proposition}

\begin{proof}
  This follows directly from the above lemma by taking
  \begin{align*}
	\cA &= \coh{X}, & \cB &= \coh{Z}, \\
	\text{ and } A^{\bullet} &= \cA_{\alpha} \text{ (concentrated at degree
	0)} &  &
  \end{align*}
  in lemma \ref{lem:tensring} where $\coh{\_}$ is the category of coherent
  sheaves over the corresponding scheme. By proposition \ref{prp:compmond},
  we get that $\moduledg{(\bM_d \bM_A)}$ and $\moduledg{(\bM_A \bM_d)}$ are
  equivalent. Part 2 of the same proposition says that these are
  equivalent to $\moduledg{\bar{\bM}_d}$ where $\bar{\bM}_d$ is the induced
  monad on $\moduledg{\bM_A}$. But $\moduledg{\bM_A}$ is equivalent to
  $\cdgr{\coh{X, \alpha}}{Z}$, of complexes whose homologies are supported
  in $Z$.

  Now $\moduledg{\bar{\bM}_d}$ is the category of all K-injective objects.
  Thus $H^0(\moduledg{(\bM_d \bM_A)})$ is nothing but $\cD_Z(X, \alpha)$.

  This proves the corollary.
\end{proof}

For the second one, we use the notations from subsection \ref{ssc:elagincst}
and section \ref{sec:bousfield}.

\tododone{Modify from Elagin: ($(T/T_Z)^G \cong T^G / T_Z^G$, support
version of equivariant der cat. Not clear if $\cD_Z$ comes from a Drinfeld
localization.}

\begin{proposition}
  Let $G$ be a finite group which acts on two triangulated categories $\cC$
  and $\cD$.  Suppose that $\cC$ and $\cD$ are $k$-linear for some
  $\Z[1/\abs{G}]$ algebra $k$. Assume that $\cC$ (resp. $\cD$) admits an
  $G$-equivariant enhancement $(\dC, \varepsilon_{\dC})$ (resp. $(\dD,
  \varepsilon_{\dD})$. Further assume that there is a $G$-equivariant
  functor $\bL : \dC \to \dD$ which is also a Drinfeld localization functor.
  Then $\cD^G$ admits an enhancement $\bL (\tilde{Q}(B))$.
\end{proposition}
\begin{proof}
   We recall from the proof of theorem \ref{thm:geqvtcat} that $\cC^G$
   admits an enhancement $(\tilde{Q}(B), \Gamma(B))$ where $B = (\dC,
   \perf{(\dC^G)}, \boldp_*^{\dC}, \boldp^*_{\dC})$. Since $\bL$ is
   $G$-equivariant, it induces a functor $\perf{(\bL^G)} : \perf{(\dC^G)}
   \to \perf{(\dC^G)}$ and hence a morphism $\bar{\bL} = (\bL,
   \perf{(\bL^G)}) : B \to B$ in $\cB$. Note $B' = \bar{\bL}(B) = (\dD,
   \perf{(\dD^G)}, \boldp_*^{\dD}, \boldp^*_{\dD})$. It is clear that $\bL
   \circ \tilde{Q} = \tilde{Q} \circ \bar{\bL}$, by construction. Thus, on
   one hand we have
   \begin{equation*}
     H^0(\bL \circ \tilde{Q} (B)) \cong
     H^0(\bL) \circ H^0(\tilde{Q}(B)) \cong L (\cC^G).
   \end{equation*}
   and on the other,
   \begin{equation*}
     H^0 (\tilde{Q} (\bar{\bL}(B))) \cong  H^0(\tilde{Q} (B')) \cong \cD^G.
   \end{equation*}
   This completes the proof.
\end{proof}

\bibliographystyle{amsalpha}
\bibliography{references}

\end{document}